\documentclass[a4paper,11pt]{article}

\usepackage[T2A]{fontenc}
\usepackage[cp1251]{inputenc}
\usepackage[english,russian]{babel}
\usepackage[tbtags]{amsmath}
\usepackage{amsfonts,amssymb,mathrsfs,amscd}
\usepackage{amsthm}
\usepackage{euscript,textcomp,verbatim,fancyhdr}
\usepackage[all]{xy}
\usepackage{epsfig}
\usepackage{latexsym}
\usepackage{mdwtab}           
\usepackage{graphics}

\numberwithin{equation}{section}
\theoremstyle{plain}
\newtheorem{Thm}{Теорема}[section]
\newtheorem{Lem}[Thm]{Лемма}
\newtheorem{Cor}[Thm]{Следствие}
\theoremstyle{definition}

                             \newtheorem{Rem}[Thm]{Замечание}
                             
                             \newtheorem{Def}[Thm]{Определение}
                             \newtheorem{Not}[Thm]{Обозначение}

\def\Im{\operatorname{Im}}
\def\const{\operatorname{const}}
\def\RR{\mathbb R}

\def\NN{\mathbb N}

\newcommand{\lam}{\lambda}
\newcommand{\ZZ}{\mathbb Z}
\newcommand{\EE}{\mathbb{E}}
\newcommand{\FF}{\mathbb{F}}
\newcommand{\XX}{\mathbb{\widetilde X}}
\newcommand{\YY}{\mathbb{\widetilde Y}}
\newcommand{\uups}{\boldsymbol{\ups}}
\newcommand{\barc}{{\boldsymbol{c}}}

\renewcommand{\SS}{\mathbb S}
\newcommand{\DD}{\mathbb D}
\newcommand{\QQ}{\mathbb{Q}}
\newcommand{\UU}{\mathbb{U}}
\newcommand{\KK}{\mathbb{\widetilde K}}
\newcommand{\IKK}{\mathbb{K}}
\newcommand{\T}{{\cal T}}
\renewcommand{\tilde}{\widetilde}
\renewcommand{\hat}{\widehat}
\newcommand{\tildeGamma}{\Theta}      
\newcommand{\K}{\widetilde K}
\renewcommand{\P}{{\cal P}}
\newcommand{\MM}{{\widetilde{\cal M}}}
\newcommand{\D}{{\mathscr D}}
\newcommand{\N}{{\cal C}}
\newcommand{\Diff}{{\rm Diff}}

\newcommand{\Homeo}{{\rm Homeo}}
\newcommand{\Hom}{{\rm Hom}}
\newcommand{\Aut}{{\rm Aut}}
\renewcommand{\Im}{{\rm Im}}
\newcommand{\num}{{\rm num}}
\newcommand{\fix}{{\rm fix}}
\newcommand{\adm}{{\rm adm}}

\newcommand{\diskr}{{\rm discr}}

\newcommand{\ind}{{\rm ind\,}}
\newcommand{\eps}{\varepsilon}
\newcommand{\ups}{\upsilon}
\newcommand{\Int}{{\rm int\,}}
\newcommand{\grad}{{\rm grad\,}}
\newcommand{\stab}{{\rm stab}}
\newcommand{\st}{{\rm st}}
\newcommand{\rank}{{\rm rank\,}}
\newcommand{\sgn}{{\rm sgn\,}}

\newcommand{\id}{{\rm id}}
\newcommand{\isot}{{\rm isot}}
\newcommand{\glue}{{\rm glue}}
\newcommand{\eval}{{\rm Ev}}
\renewcommand{\mod}{{\rm \,mod\,}}

\renewcommand{\:}{\colon\,}

\renewcommand{\b}{y}

\renewcommand{\d}{\partial}

\renewcommand{\o}[1]{{\stackrel{{\scriptscriptstyle\circ}}{#1}}}
\newcommand{\<}{\langle}
\renewcommand{\>}{\rangle}
\newcommand{\aapprox}{\sim} 
\newcommand{\ccong}{\approx} 
\renewcommand{\emptyset}{\varnothing}
\begin{document}

\title{On the homotopy type of the spaces of Morse functions on surfaces
}
\author
{E.\,A.~Kudryavtseva}
\date{}

\maketitle


Let $M$ be a smooth closed orientable surface. Let $F$ be the space
of Morse functions on $M$ having fixed number of critical points of
each index, moreover at least $\chi(M)+1$ critical points are labeled
by different labels (enumerated). A notion of a skew
cylindric-polyhedral complex, which generalizes the notion of a
polyhedral complex, is introduced. The skew cylindric-polyhedral
complex $\KK$ (the ``complex of framed Morse functions''), associated
with the space $F$, is defined. In the case when $M=S^2$, the
polyhedron $\KK$ is finite; its Euler characteristic $\chi(\KK)$ is
evaluated and the Morse inequalities for its Betti numbers
$\beta_j(\KK)$ are obtained. A relation between the homotopy types of
the polyhedron $\KK$ and the space $F$ of Morse functions, endowed
with the $C^\infty$-topology, is indicated.



{\bf MSC-class:} 58E05, 57M50, 58K65, 46M18

\markright{On the homotopy type of the spaces of Morse functions}

\footnotetext[0]{Работа выполнена при поддержке РФФИ (грант
\No~10–01–00748-а), Программы поддержки ведущих научных школ РФ
(грант \No~НШ-3224.2010.1), Программы ``Развитие научного потенциала
высшей школы'' (грант \No~2.1.1.3704), ФЦП ``Научные и
научно-педагогические кадры инновационной России'' (грант
\No~14.740.11.0794).}

УДК 515.164.174+515.164.22+515.122.55

\begin{center}
{\LARGE О гомотопическом типе пространств функций Морса \\ на поверхностях} \bigskip \\
\large Е.\,А.~Кудрявцева
\bigskip
\end{center}

\begin{abstract}
Пусть $M$ --- гладкая замкнутая ориентируемая поверхность. Пусть $F$
-- пространство функций Морса на $M$ с фиксированным количеством
критических точек каждого индекса, причем не менее чем $\chi(M)+1$
критических точек помечены различными метками (пронумерованы).
Введено понятие косого цилиндрически-полиэдрального комплекса,
обобщающее понятие полиэдрального комплекса.
Определен косой цилиндрически-поли\-эд\-раль\-ный комплекс $\KK$
(``комплекс оснащенных функций Морса''), ассоциированный с
пространством $F$.
В случае $M=S^2$ полиэдр $\KK$ конечен; вычислена его эйлерова
характеристика $\chi(\KK)$ и получены неравенства Морса для его чисел
Бетти $\beta_j(\KK)$. Указана связь гомотопических типов полиэдра
$\KK$ и пространства $F$ функций Морса, снабженного
$C^\infty$-топологией.

Библиография: 41 название.

{\bf Ключевые слова:} функции Морса, комплекс оснащенных функций Морса, полиэдральный комплекс,
$C^\infty$-топология, универсальное пространство модулей.
\end{abstract}

\section{Введение} \label{sec:intro}

Настоящая работа является продолжением работы~\cite {kp1}.

Задача изучения гладких функций с ``умеренными'' особенностями на
гладком многообразии $M$ является классической. Изучение топологии
пространства таких функций как правило состоит из двух частей:
 \begin{itemize}
 \item [1)] сведение к комбинаторной задаче, т.е.\ построение комбинаторного
объекта (например, конечномерного полиэдра), гомотопически
эквивалентного изучаемому пространству функций;
 \item[2)] решение полученной комбинаторной задачи (изучение топологии
построенного полиэдра).
 \end{itemize}
Одним из таких подходов является (параметрический) $h$-принцип
М.~Громова~\cite{Gromov86}, изучаемый, например, в работах
\cite{Igusa84}, \cite{Vassil89}, \cite[теорема 2.3]{Franks},
\cite{Hatcher75,ChL09} (см.\ также \cite {kp1}).
Невыполнение 1-параметрического $h$-принципа для пространств функций Морса на некоторых компактных многообразиях размерности большей $5$ показано в работах~\cite {ChL09,Hatcher75} (см.\ также \cite [\S1]{kp1}).

Рассмотрим задачу о вычислении гомотопического типа пространства
$F(M)$ функций Морса на компактном гладком многообразии $M$,
например, на гладкой поверхности. Для пространства $F(S^1)$ функций
Морса на окружности $M=S^1$ имеется следующий метод решения.
Сопоставим каждой функции Морса $f\in F_r(S^1)$, имеющей ровно $2r$
критических точек,
множество ее критических точек локальных минимумов (т.е.\ некоторую
$r$-точечную конфигурацию на поверхности $M$). Нетрудно доказывается,
что построенное отображение $F_r(S^1)\to Q_r(S^1)$ сюръективно и
является гомотопической эквивалентностью. Тем самым, описанный метод
сводит задачу к изучению топологии конфигурационного пространства
$Q_r(S^1)$, т.е.\ пространства $r$-точечных конфигураций на $S^1$.
Гомотопический тип последнего пространства легко находится и равен
$S^1$.

В настоящей работе изучается топология пространства $F=F(M)$ функций
Морса на компактной двумерной поверхности $M$. Предлагаемый нами
метод аналогичен ``методу конфигурационных пространств'', описанному
выше. А именно, в настоящей работе
описывается построение комбинаторного объекта -- комплекса
$\KK=\KK(M)$ оснащенных функций Морса (определение~\ref {def:KK} и
теорема~\ref {thm:KK}), аналогичного комплексу функций Морса
$\K=\K(M)$ из работы~\cite {BaK}.
Комплекс $\KK$ является конечномерным косым
цилиндрически-полиэдральным комплексом (см.\ определение~\ref
{def:pol}), т.е.\ допускает разбиение на ``косые цилиндрические
ручки'' (см.\ определение~\ref {def:thick:cylinder}), аналогичные
круглым ручкам, и приклеенные друг к другу ``регулярным'' образом.
При этом ручки находятся во взаимно однозначном соответствии с
классами изотопности
функций Морса из $F^1\subset F$ (см.\ определения~\ref
{def:Morse}(B), \ref {def:equiv}), а подошва ручки $\DD_{[f]_\isot}$,
отвечающей классу изотопности $[f]_\isot$ функции $f$, содержится в
объединении ручек $\DD_{[g]_\isot}$, отвечающих классам изотопности
функций, полученных малыми возмущениями функции $f$.
Важным свойством комплекса $\KK$ является то, что
в большинстве случаев (см.\ (\ref {eq:main}))
пространство $F$ функций Морса на компактной поверхности $M$
гомотопически эквивалентно 
полиэдру $R\times\KK$: 
 \begin{equation} \label {eq:sim}
F\sim R\times\KK,
 \end{equation}
где $R=R(M)$ -- одно из многообразий $\RR P^3$, $S^1$, $S^1\times
S^1$ и точка (см.\ (\ref {eq:EE})), а $\KK=\KK(M)$ -- ``комплекс
оснащенных функций Морса'',
построенный в настоящей работе (см.\ замечание \ref {rem:thm}).
Тем самым,
изучение топологии пространства $F$ функций Морса сводится к
комбинаторной задаче --- изучению топологии 
полиэдра $\KK$.

В некоторых случаях гомологии полиэдра $\KK$ могут быть изучены
методами теории Морса, ввиду естественного разложения полиэдра $\KK$
на косые цилиндрические ручки.
В качестве иллюстрации мы получаем обобщенные неравенства Морса для
чисел Бетти пространства
$\KK$ и находим его эйлерову характеристику в случае, когда $M=S^2$ и
у каждой функции $f\in F$ не менее трех критических точек помечены
разными метками, т.е.\ занумерованы (следствие \ref {cor:ineq}).

Наши комплексы $\K$ и $\KK$ функций Морса и оснащенных функций Морса
аналогичны известным (абстрактным симплициальным флаговым)
комплексам, рассматриваемым при изучении группы классов отображений
поверхности $M$ (см.~\cite {FI}, \cite{H,Ha,I1}, \cite{HT}, \cite
{Mar}), кубическим комплексам~\cite{MP}, и обобщают граф разрезов
(см.~\cite {Kmsb}). В случае $M=S^1$ аналогичный комплекс $\KK$
состоит из одной точки и $R=S^1$.

Связные компоненты пространств функций Морса на поверхности изучались
С.В.\ Матвеевым~\cite{Kmsb}, Х.~Цишангом, В.В.\ Шарко~\cite{SH},
Е.А.\ Кудрявцевой~\cite{Kmsb}, С.И.\ Максименко~\cite{Max2005}, а
также Ю.М.\ Бурманом~\cite{Bu,Bu2} (для пространств гладких функций
без критических точек на открытых поверхностях) и Е.А.\
Кудрявцевой~\cite{BaK} (для пространств функций Морса с фиксированным
множеством критических точек). Комбинаторные и топологические
свойства пространств функций Морса на поверхностях изучались в
работах~\cite {Kul}, \cite{Max}. В работах \cite {F86}, \cite{FZ0},
\cite{BF0,BFbook}, \cite{K}, \cite{K:stab} и \cite{BrK:stab} функции
Морса изучались в связи с задачей классификации (лиувиллевой,
орбитальной) не\-вы\-рожд\-ен\-ных интегрируемых гамильтоновых систем
с двумя степенями свободы. Группы гомологий и гомотопий пространств
функций с умеренными особенностями (с допущением не-морсовских
особенностей) на окружности изучался В.И.\ Арнольдом~\cite{A89}.

Перейдем к точным формулировкам.

 \begin{Def} \label {def:Morse}
Пусть $M$ --- гладкая (т.е.\ класса $C^\infty$) компактная связная
поверхность, край которой пуст или не пуст, $\d M=\d^+M\sqcup\d^-M$,
где $\d^+M$ --- объединение некоторых граничных окружностей $M$.

{\rm(A)} Обозначим через $C^\infty(M)$ пространство гладких (т.е.\
класса $C^\infty$) вещественнозначных функций $f$ на $M$. Обозначим
через $C^\infty(M,\d^+M,\d^-M)\subset C^\infty(M)$ подпространство,
состоящее из таких функций $f\in C^\infty(M)$, что все ее критические
точки (т.е.\ такие точки $x\in M$, что $df|_x=0$) принадлежат $\Int
M$, а любая граничная точка $x\in\d M$ имеет такую окрестность $U$ в
$M$, что $f(U\cap\d M)=f(x)$, причем $\inf(f|_U)=f(x)$ при
$x\in\d^-M$, и $\sup(f|_U)=f(x)$ при $x\in\d^+M$.

{\rm(B)} Функцию $f$ на $M$ назовем {\em функцией Морса на
$(M,\d^+M,\d^-M)$}, если $f\in C^\infty(M,\d^+M,\d^-M)$ и все ее
критические точки невырождены (т.е.\ квадратичная форма в $T_xM$,
заданная матрицей вторых частных производных $f$ в критической точке
$x$, невырождена). Пусть
$$
F:=F_{p,q,r}(M,\d^+M,\d^-M)
$$
--- пространство функций Морса $f$ на поверхности $(M,\d^+M,\d^-M)$,
имеющих ровно $p+q+r$ критических точек, из которых $p$ точек
локальных минимумов, $q$ седловых точек и $r$ точек локальных
максимумов. Пусть $d^+,d^-\ge0$ --- число граничных окружностей в
$\d^+M$ и $\d^-M$ соответственно. Будем предполагать, что выполнены
неравенства Морса:
 \begin{equation} \label {eq:Morse:ineq}
\chi(M)=p-q+r, \quad p+d^+>0, \quad r+d^->0,
 \end{equation}
так как в противном случае $F=\emptyset$. Пространство $F$ мы наделим
$C^\infty$-топо\-ло\-ги\-ей, см.~\cite[\S4(а)]{kp1}, и назовем его {\em
пространством функций Морса на поверхности $(M,\d^+M,\d^-M)$}.
Обозначим через $F^1\subset F$ подпространство в $F$, состоящее из
таких функций Морса $f\in F$, что все локальные минимумы равны
$f(\d^-M)=-1$, а все локальные максимумы равны $f(\d^+M)=1$.

{\rm(C)} Обозначим через $F^\num$ (соотв.\ $F^{1,\num}$)
пространство, полученное из пространства $F$ (соотв.\ $F^1$)
введением нумерации у некоторых из критических точек (называемых {\em
отмеченными} критическими точками) функций Морса $f\in F$
(соответственно $f\in F^1$). Наделим его $C^\infty$-топологией как
в~\cite[\S4(а)]{kp1}. Если количество отмеченных критических точек
локальных минимумов, максимумов и седловых точек равно $\hat p,\hat
r,\hat q$ соответственно (где $0\le\hat p\le p$, $0\le\hat q\le q$,
$0\le\hat r\le r$), то имеем $C_p^{\hat p}C_q^{\hat q}C_r^{\hat
r}\hat p!\hat q!\hat r!$-листные накрытия $F^\num\to F$ и
$F^{1,\num}\to F^1$.
 \end{Def}

\begin{Not} \label{not:R:G0*}
Пусть $\D^\pm=\Diff(M,\d^+M,\d^-M)$ --- группа всех (не обязательно сохраняющих
ориентацию и компоненты края) диффеоморфизмов поверхности $M$, переводящих каждое множество $\d^+M,\d^-M$ в себя.
Пространство $\D^\pm$ наделим $C^\infty$-топологией,
см.~\cite[\S4(б)]{kp1}. Пусть $\D^0\subset\D^\pm$ --- подгруппа,
состоящая из всех диффеоморфизмов $h\in\D^\pm$, гомотопных $\id_M$ в
классе гомеоморфизмов $M$.
 \end{Not}

\begin{Def} \label{def:equiv}
{\rm(А)} Две функции Морса $f,g\in F$ назовем {\em эквивалентными},
если найдутся такие диффеоморфизмы $h_1\in\D^\pm$ и
$h_2\in\Diff^+(\RR)$, что $f=h_2\circ g\circ h_1$
 (и $h_1$ сохраняет нумерацию критических точек, если $f,g\in F^{\num}$),
и обозначаем $f\sim g$. Класс эквивалентности функции $f$ обозначим
через $[f]$.

{\rm(Б)} Две функции Морса $f$ и $g$ назовем {\em изотопными}, если
они эквивалентны и $h_1\in\D^0$ (т.е.\ $h_1$ изотопен
тождественному), и обозначаем $f\sim_\isot g$. Множество всех функций
из $F^1$, изотопных функции $f$, обозначим через $[f]_\isot$.
\end{Def}

Классификация функций Морса из $F$ с точностью до {\em
эквивалентности} легко получается из классификации функций Морса с
точностью до {\em послойной эквивалентности} (см.~\cite[гл.\,2,
теоремы~4 и~8]{BFbook}). Автором доказаны критерии (послойной)
эквивалентности (соотв.\ изотопности) пары ``возмущенных'' функций
Морса $\tilde f,\tilde g\in C^\infty(M,\d^+M,\d^-M)$, полученных
малыми возмущениями из пары (послойно) эквивалентных (соотв.\
изотопных) функций Морса $f,g\in F$ (см.\ \S\ref {sec:KK0}, шаг 3,
или \cite[утверждение~1.1 и~\S3]{K}), и доказана классификация
функций Морса с точностью до {\it изотопности} (см.~\cite[лемма~1 и
теорема~2]{KP2}).

В работе~\cite{kp1} введено понятие {\it оснащенных} функций Морса на
компактной поверхности $M$
и доказана гомотопическая эквивалентность пространства $F$
функций Морса и пространства $\FF$ оснащенных функций Морса на $M$.
Также доказан аналог последнего утверждения для
ограничений указанных гомотопических эквивалентностей на классы
изотопности $[f]_\isot$ (\cite[теорема 2.5]{kp1}). Последнее дает
положительный ответ на вопрос, поставленный В.И.~Арнольдом.

Статья имеет следующую структуру. В~\S \ref {subsec:analyt}
определяется более общее пространство функций Морса, вводится понятие
косого цилиндрически-поли\-эд\-раль\-ного комплекса и формулируются
основные результаты настоящей работы (теорема~\ref {thm:KP4} и
следствие \ref {cor:ineq}).
В \S\ref {sec:KK0} описывается построение стандартной косой
цилиндрической ручки $\DD_{[f]_\isot}^\st$, отвечающей классу
изотопности $[f]_\isot$ функции Морса $f\in F^1$, и изучены
отображения инцидентности между парами ручек, отвечающими парам
примыкающих друг к другу классов изотопности
$[f]_\isot\prec[g]_\isot$ функций (определение~\ref {def:pol}(D),
леммы \ref{lem:P}--\ref {lem:tilde}).
В \S\ref {sec:KK} описывается построение комбинаторного объекта --
комплекса $\KK$ оснащенных функций Морса
(и содержащего его многообразия $\MM$) и доказывается, что он
является косым цилиндрически-полиэдральным комплексом (теоремы \ref
{thm:KK}, \ref {thm:MM1}).

Автор приносит благодарность В.И.~Арнольду, С.А.\ Мелихову,
А.А.~Ошемкову, Д.А.~Пермякову, А.Т.~Фоменко и Х.~Цишангу за внимание к работе и
полезные обсуждения.

\section{Точные формулировки основных результатов}\label{subsec:analyt}

Следующее пространство обобщает пространства
$F_{p,q,r}(M,\d^+M,\d^-M)$, $F^\num$ и состоит из функций Морса, у
которых некоторые из отмеченных критических точек закреплены.

\begin{Def} \label{def:F}
{\rm(A)} Пусть $0\le p^*\le\hat p$, $0\le q^*\le\hat q$, $0\le r^*\le\hat r$ (см.\
определение \ref {def:Morse}(C)). Обозначим
 $$
(p',p'';q',q'';r',r''):=(\hat p-p^*,p-\hat p;\hat q-q^*,q-\hat q;\hat
r-r^*,r-\hat r).
 $$
Для каждой функции $f\in F^\num$ обозначим через $\N_{f,0}$,
$\N_{f,1}$, $\N_{f,2}$ множества ее критических точек локальных
минимумов, седловых критических точек и точек локальных максимумов
соответственно, и через $\hat\N_{f,\lam}\subseteq \N_{f,\lam}$,
$\lam=0,1,2$, множества отмеченных критических точек. В каждом из
множеств отмеченных (а потому занумерованных) критических точек
рассмотрим подмножество, обозначаемое через $\N^*_{f,0}$,
$\N^*_{f,1}$, $\N^*_{f,2}$ и состоящее из первых $p^*,q^*,r^*$ точек
соответственно. Фиксируем ``базисную'' функцию $f_*\in F^\num$. Пусть
$$
F:=F_{p^*,p',p'';q^*,q',q'';r^*,r',r''}(M,\d^+M,\d^-M)
$$
-- пространство функций Морса $f\in F^\num$ на поверхности
$(M,\d^+M,\d^-M)$ (см.\ определение~$\ref {def:Morse}${\rm(B)}),
таких что $\N^*_{f,\lambda}=\N^*_{f_*,\lambda}$ для любого
$\lambda=0,1,2$. Пространство $F$ мы наделим $C^\infty$-топологией,
см.~\cite[\S4]{kp1}, и назовем его {\em обобщенным пространством
функций Морса на поверхности $(M,\d^+M,\d^-M)$}. Определим
подпространство $F^1\subset F$ аналогично определению \ref {def:Morse}(B).
 \end{Def}

Пространства $F,F^\num$ из определения \ref {def:Morse} имеют вид
$$
F =F_{0,0,p;0,0,q;0,0,r}(M,\d^+M,\d^-M), \quad F^\num=F_{0,\hat
p,p'';0,\hat q,q'';0,\hat r,r''}(M,\d^+M,\d^-M)
$$
соответственно. Из теоремы С.В.\ Матвеева (см.~\cite{Kmsb}) и ее
обобщения в~\cite {Kmsb} для пространства $F^\num$ получаем следующее
ее обобщение: любое обобщенное пространство $F=F_{p^*,p',p'';0,\hat
q,q'';r^*,r',r''}(M,\d^+M,\d^-M)$ функций Морса без закрепленных
седловых точек (т.е.\ при $q^*=0$) линейно связно.

Всюду в статье мы предполагаем, что поверхность $M$ ориентирована.
Случай неориентируемой поверхности $M$ рассматривается как в
\cite[замечание 2.7]{kp1}.

Следующие группы диффеоморфизмов обобщают группы $\D^\pm$ и $\D^0$,
см.\ обо\-зна\-че\-ние~\ref {not:R:G0*} выше. Мы их вводим для
изучения обобщенного пространства $F$ функций Морса (и по-прежнему
обозначаем их через $\D^\pm$ и $\D^0$).

\begin{Not} \label{not:R:G0}
(A) Множества $\N^*_{f,0},\N^*_{f,1},\N^*_{f,2}$ фиксированных
критических точек совпадают для разных функций $f\in F$, будем их
обозначать через $\N_0,\N_1,\N_2$ соответственно, положим
$\N:=\N_0\cup\N_1\cup\N_2$. Обозначим через $\D^\pm=\Diff(M,\d^+M,\d^-M,\N_0,\N_1,\N_2)$
группу всех (не обязательно сохраняющих ориентацию и компоненты
края) диффеоморфизмов поверхности $M$, переводящих каждое множество $\d^+M,\d^-M$,
$\N_\lam$ в себя, $0\le\lam\le2$. Пространство $\D^\pm$ тоже наделим
$C^\infty$-топологией, см.~\cite[\S4(б)]{kp1}. Если $M$ ориентируема,
обозначим через $\D\subset\D^\pm$ подгруппу сохраняющих ориентацию
диффеоморфизмов (с индуцированной топологией). Пусть
$\D^0=\Diff^0(M,\N)\subset\D^\pm$ --- подгруппа (с
индуцированной топологией), состоящая из всех диффеоморфизмов
$h\in\D^\pm$, гомотопных $\id_M$ в классе гомеоморфизмов
пары $(M,\N)$.

(B) Обозначим через $\bar M$ замкнутую поверхность, полученную из
поверхности $M$ стягиванием в точку каждой граничной окружности.
Обозначим через $\T\subset\D$ группу (называемую {\em группой
Торелли}), состоящую из всех диффеоморфизмов $h\in\D$, переводящих в
себя каждую компоненту края $M$, и таких что индуцированный
гомеоморфизм $\bar h\:\bar M\to\bar M$ индуцирует тождественный
автоморфизм группы гомологий $H_1(\bar M)$. Имеем $\D^0\subset\T$.
 \end{Not}

Из результатов~\cite {EE,EE0} К.Дж.~Эрля и Дж.~Иллса (мл.) следует,
что имеется гомотопическая эквивалентность
 \begin{equation} \label {eq:EE}
 \D^0 \aapprox R_{\D^0},
 \end{equation}
где $R_{\D^0}$ --- одно из четырех многообразий, определяемое парой
$(M,|\N|)$, а именно: $SO(3)=\RR P^3$ (при
$M=S^2$, $\N=\emptyset$), $SO(2)=S^1$ (при $0\le\chi(M)-|\N|\le1$ и
$d^++d^-+|\N|>0$), $T^2=S^1\times S^1$ (при
$M=T^2$, $\N=\emptyset$) и точка (при $\chi(M)<|\N|$) (см.,
например,~\cite{S,EE}). В частности, группа $\D^0$ линейно связна.

\subsection{Косые цилиндрически-полиэдральные комплексы}
\label{subsec:KK}

Всюду в статье многогранники выпуклы и имеют размерность $\le2q$, а
евклидовы многогранники $\le q-1$. Под {\it выпуклым многогранником}
(соответственно {\it евклидовым выпуклым многогранником}) понимаем
выпуклую оболочку конечного подмножества векторного пространства
$\RR^{2q}$ (соответственно евклидова пространства $\EE^{q-1}$), а под
{\it изоморфизмом} многогранников --- биекцию между многогранниками,
продолжающуюся до аффинного изоморфизма (соответственно изометрии)
объемлющих пространств.

{\it Утолщенный цилиндр} --- это главное расслоение над выпуклым
многогранником со слоем стандартный цилиндр $\RR^c\times(S^1)^d$ (где
прямые сомножители $S^1$ в разложении цилиндра не упорядочены), а
{\it стандартная цилиндрическая ручка} --- это прямое произведение
евклидова выпуклого многогранника и утолщенного цилиндра. Уточним и
обобщим эти понятия.

\begin{Def} [стандартная косая цилиндрическая ручка] \label {def:thick:cylinder}
(A) {\em Утолщенным цилиндром} назовем прямое произведение
$\SS:=(\RR^c\times(S^1)^d)\times P$, где $P$ --- выпуклый
многогранник, $S^1=\RR/\ZZ$, $c,d\in\NN\cup\{0\}$. Гомеоморфизм
$h:\SS_1\to\SS_2$ между утолщенными цилиндрами назовем {\em
допустимым}, если $c_1=c_2=:c$, $d_1=d_2=:d$, существуют биекции
$\pi\in\Sigma_c$ и $\rho\in\Sigma_d$, изоморфизм многогранников
$a\:P_1\to P_2$ и непрерывное отображение
$\delta=(\delta_1,\dots,\delta_{c+d})\:P_1\to\RR^c\times(S^1)^d
$, такие что для любого $(x_1,\dots,x_c,\varphi_1,\dots,\varphi_d,u)\in\SS_1$ выполнено
 $$
 h(x_1,\dots,x_c,\varphi_1,\dots,\varphi_d,u)
 =(x_{\pi(1)},\dots,x_{\pi(c)},\varphi_{\rho(1)},\dots,\varphi_{\rho(d)},0)
 $$
 $$
 +(\delta_{1}(u),\dots,\delta_{c+d}(u),a(u)).
 $$

Автоморфизм $b\:D\to D$ евклидова многогранника (т.е.\ изоморфизм на
себя) назовем {\it допустимым}, если он тривиален или не имеет
неподвижных вершин, а для любой грани $\tau\subset\d D$ выполнено
либо $b(\tau)=\tau$, либо $\tau\cap b(\tau)=\emptyset$.

(B) {\em Стандартной цилиндрической ручкой} назовем прямое
произведение $D\times\SS$ евклидова выпуклого многогранника $D$ и
утолщенного цилиндра $\SS$ (см.~(А)). Гомеоморфизм $D_1\times\SS_1\to
D_2\times\SS_2$ между стандартными цилиндрическими ручками назовем
{\em изоморфизмом}, если он является прямым произведением изоморфизма
$b\:D_1\to D_2$ евклидовых многогранников и допустимого гомеоморфизма
$\SS_1\to\SS_2$ утолщенных цилиндров (см.\ (A)).

Автоморфизм $D\times\SS\to D\times\SS$ стандартной цилиндрической
ручки назовем {\it допустимым}, если либо он совпадает с тождественным,
либо хотя бы один из соответствующих автоморфизмов многогранников $b\:D\to
D$, $a\:P\to P$ и перестановок $\pi\in\Sigma_c$ и $\rho\in\Sigma_d$
нетривиален,
причем автоморфизм $b\:D\to D$ допустим (см.\ (A)) и выполнены
следующие (необязательные в общем случае, но выполненные для
комплексов оснащенных функций Морса в случае (\ref {eq:main}))
дополнительные условия: перестановка $\pi\in\Sigma_c$ всегда
тривиальна, в случае тривиальности автоморфизма $a\:P\to P$
автоморфизм $b\:D\to D$ тривиален, а в случае тривиальности $b$
перестановка $\rho\in\Sigma_d$ и автоморфизм $a^2\:P\to P$
тривиальны.

(C) {\em Стандартной косой цилиндрической ручкой}
назовем пространство орбит $\DD^\st:=(D\times\SS)/\Gamma$ свободного
действия (автоматически конечной) группы $\Gamma$ на стандартной
цилиндрической ручке $D\times\SS$ допустимыми автоморфизмами
(см.~(B)). Размерность $k$ многогранника $D=D^k$ назовем {\em
индексом}
ручки $\DD=\DD^\st$, подмножество $\d\DD:=((\d
D)\times\SS)/\Gamma\subset\DD$ назовем ее {\em подошвой}, а
дополнение $\o\DD:=\DD\setminus\d\DD$ --- {\em открытой стандартной
(косой) цилиндрической ручкой}, отвечающей
ручке $\DD$. Для любой грани $D'\subset\d D$ образ подмножества
$D'\times\SS\subset D\times\SS$ при проекции $D\times\SS\to\DD^\st$
назовем {\em (косой) гранью} стандартной (косой) ручки $\DD^\st$.
(Косая грань всегда является стандартной косой цилиндрической ручкой
ввиду допустимости автоморфизма $b\:D\to D$ из (B), см.\ (A).)
Гомеоморфизм $\DD_{1}^\st\to\DD_{2}^\st$ между стандартными косыми
цилиндрическими ручками назовем {\em изоморфизмом}, если он
поднимается до изоморфизма $D_1\times\SS_1\to D_2\times\SS_2$
соответствующих стандартных цилиндрических ручек. Автоморфизм
$\DD^\st\to\DD^\st$ стандартной косой цилиндрической ручки назовем
{\it допустимым}, если он поднимается до допустимого автоморфизма
$D\times\SS\to D\times\SS$ соответствующей стандартной цилиндрической
ручки (см.\ (B)).

(D) Погружение $i:\SS_1\looparrowright\SS_2$ между утолщенными
цилиндрами $\SS_j=(\RR^{c_j}\times(S^1)^{d_j})\times P_j$, $j=1,2$
(см.~(A)), назовем {\em допустимым}, если существуют отображения
$\pi\:\{1,\dots,c_2+d_2\}\to\{0,1,\dots,c_1+d_1\}$ и
$\eta\:\{1,\dots,c_1+d_1\}\to\{1,-1\}$,
такие что
$\{1,\dots,c_1\}\subset\pi(\{1,\dots,c_2\})\subset\{0,1,\dots,c_1\}$,
$\pi|_{\{1,\dots,c_2\}\cap\pi^{-1}\{1,\dots,c_1\}}$ инъективно,
$\{c_1+1,\dots,c_1+d_1\}\subset\pi(\{c_2+1,\dots,c_2+d_2\})$,
$\eta(\{1,\dots,c_1\})=1$, а также существуют отображение $a\:P_1\to
P_2$, продолжающееся до аффинного, и непрерывное отображение
$\delta=(\delta_1,\dots,\delta_{c_2+d_2})\:P_1\to\RR^{c_2}\times(S^1)^{d_2}
$, такие что
$$
 i(x_1,\dots,x_{c_1+d_1},u)
 =(\eta(\pi(1))\,x_{\pi(1)},\dots,\eta(\pi(c_2+d_2))\,x_{\pi(c_2+d_2)},0)
$$
$$
 +(\delta_{1}(u),\dots,\delta_{c_2+d_2}(u),a(u))
$$
для любого $(x_1,\dots,x_{c_1+d_1},u)\in\SS_1$, где последние $d_j$
координат любой точки из $\RR^{c_j}\times(S^1)^{d_j}$ рассматриваются
по модулю 1 при $j=1,2$, $x_0:=0$, $\eta(0):=1$. Погружение
$D_1\times\SS_1\looparrowright D_2\times\SS_2$ между стандартными
цилиндрическими ручками (см.\ (B)) назовем {\em допустимым}, если оно
является прямым произведением изоморфизма $D_1\to D_2$ евклидовых
многогранников и допустимого погружения $\SS_1\looparrowright\SS_2$
утолщенных цилиндров. Вложение $\DD_1^\st\hookrightarrow\DD_2^\st$
между стандартными {\it косыми} цилиндрическими ручками (см.~(C))
назовем {\em мономорфизмом}, если оно поднимается до допустимого
погружения $D_1\times\SS_1\looparrowright D_2\times\SS_2$
соответствующих стандартных цилиндрических ручек.
\end{Def}

Пусть $X$ ---
топологическое пространство.

\begin {Def} [косой цилиндрически-полиэдральный комплекс] \label {def:pol}
(A) Будем говорить, что на подмножестве $\DD\subset X$ задана {\em
структура (косой) цилиндрической ручки}, и называть это подмножество
{\it (косой) цилиндрической ручкой},
если $\DD$ замкнуто в $X$, и фиксированы стандартная (косая)
цилиндрическая ручка $\DD^\st$ с точностью до изоморфизма, и
гомеоморфизм
$\varphi_\DD\:\DD^\st\stackrel{\approx}{\longrightarrow}\DD$
(называемый {\em характеристическим отображением} (косой) ручки
$\DD$) с точностью до допустимых автоморфизмов стандартной (косой)
ручки $\DD^\st$. Подмножество $\d\DD:=\varphi_\DD(\d\DD^\st)$ назовем
{\it подошвой} (косой) ручки $\DD$. Вложение
$i\:\DD_1\hookrightarrow\DD_2$ между (косыми) цилиндрическими ручками
назовем {\em мономорфизмом}, если вложение $\varphi_{\DD_2}^{-1}\circ
i\circ\varphi_{\DD_1}$ соответствующих стандартных (косых) ручек
является мономорфизмом.

(B) Пространство $X$ назовем {\it (косым) цилиндрически-полиэдральным
комплексом}, если фиксировано разбиение $X=\bigcup_{i=1}^n\o\DD_i$,
где $n\le\infty$, на попарно непересекающиеся подмножества $\o\DD_i$,
называемые {\it открытыми (косыми) цилиндрическими ручками}
разбиения, и для каждой открытой ручки $\o\DD_i$ фиксирована
структура (косой) цилиндрической ручки на ее замыкании
$\DD_i:=\overline{\o\DD_i}$, называемом {\it (косой) цилиндрической
ручкой}
разбиения, такая что $\o\DD_i=\DD_i\setminus\d\DD_i$, причем
выполнены следующие условия:
\begin{enumerate}
\item[(c)] для любой (косой) ручки $\DD_i$
ограничение характеристического отображения
$\varphi_{\DD_i}\:\DD_i^\st\stackrel{\approx}{\longrightarrow}\DD_i$
на произвольную (косую) грань $(\DD_i^\st)'\subset\d\DD_i^\st$
соответствующей стандартной (косой) ручки $\DD_i^\st$ является
мономорфизмом $(\DD_i^\st)'\hookrightarrow\DD_j$ в некоторую (косую)
ручку $\DD_j$ (см.\ (A) и определение~\ref
{def:thick:cylinder}(C,D)), откуда подошва $\d\DD_i$ любой (косой)
ручки $\DD_i$ индекса $k$ содержится в объединении конечного числа
(косых) ручек индекса $k-1$;
\item[(w)] подмножество $Y\subset X$ замкнуто тогда и только тогда, когда
для любой (косой) ручки $\DD_i$ замкнуто пересечение $Y\cap\DD_i$.
\end{enumerate}
Максимальный индекс (косых) цилиндрических ручек (косого)
ци\-лин\-дри\-чески\-по\-ли\-э\-драль\-ного комплекса назовем {\em
рангом} этого комплекса.

(C) Если для каждой (косой) цилиндрической ручки $\DD_i\subset X$
выполнено $c=d=\dim P=0$, получаем определение {\em строго
полиэдрального комплекса} (см.~\cite{BaK}). Определение {\em
полиэдрального комплекса} имеется, например, в~\cite{BrHae}.

(D) Пусть $\sigma,\tau\subset X$ -- два непересекающихся подмножества
топологического пространства $X$ (например, две открытые клетки
клеточного комплекса). Будем говорить, что $\sigma$ {\it примыкает к}
$\tau$ и писать $\tau\prec\sigma$ (и $\bar\tau\prec\bar\sigma$), если
$\tau\subset\d\sigma:=\bar\sigma\setminus\sigma$. Пишем
$\tau\preceq\sigma$, если $\tau\prec\sigma$ или $\tau=\sigma$.
 \end {Def}

\begin{Not} \label{not:Gf:or}
Для любой функции Морса $f\in F$ рассмотрим граф $G_f$ в поверхности
$\Int(M)$, полученный из графа $f^{-1}(f(\N_{f,1}))$ выкидыванием
всех компонент связности, не содержащих седловых критических точек
(см.\ определение \ref {def:F}). Этот граф имеет $q$ вершин
(являющихся седловыми точками $\b\in\N_{f,1}$), степени всех вершин
равны $4$, а значит в графе $2q$ ребер. Если поверхность $M$
ориентирована, то на ребрах графа $G_f$ имеется естественная
ориентация, такая, что в любой (внутренней) точке ребра репер,
составленный из положительно ориентированного касательного вектора к
ребру и вектора $\grad f$ (по отношению к какой-нибудь фиксированной
римановой метрике), задает положительную ориентацию поверхности.
Аналогично вводится ориентация на любой связной компоненте линии
уровня $f^{-1}(a)$ функции $f$, не содержащей критическую точку,
$a\in\RR$.
\end{Not}

Сформулируем основной результат данной работы, описывающий
комбинаторный объект -- комплекс $\KK$ оснащенных функций Морса,
ассоциированный с пространством $F$.

\begin{Thm}  \label{thm:KP4}
Пусть $M$ --- связная компактная ориентируемая поверхность с
разбиением края $\d M=\d^+M\sqcup\d^-M$ на положительные и
отрицательные граничные окружности. Рассмотрим обобщенные
пространства $F=F_{p^*,p',p'';q^*,q',q'';r^*,r',r''}(M,\d^+M,\d^-M)$
и $F^1\subset F$ функций Морса на поверхности $(M,\d^+M,\d^-M)$, у
которых могут быть отмеченные критические точки, из которых некоторые
точки могут быть закрепленными {\rm(см.\ определение~$\ref
{def:F}$)}.
Предположим, что
 \begin{equation} \label {eq:main}
\hat p+\hat q+\hat r >\chi(M)
 \end{equation}
(т.е.\ количество отмеченных критических точек превосходит
$\chi(M)$).
Тогда:

{\rm (A)} Имеется
косой цилиндрически-полиэдральный комплекс
 $$
 \KK=\KK_{p^*+d^-,p',p'';q^*,q',q'';r^*+d^+,r',r''}
 $$
(называемый {\em комплексом оснащенных функций Морса})
ранга $q-1$ и размерности $\dim\KK=3q-2$ при $q\ge2$ и $\dim\KK=0$
при $q\le1$, косые цилиндрические ручки которого находятся во взаимно
однозначном соответствии с классами изотопности $[f]_\isot$ функций
Морса $f\in F^1$. Индекс ручки $\DD_{[f]_\isot}$, отвечающей классу
изотопности $[f]_\isot$, равен $q-s(f)$, где $s(f)$ --- количество
седловых критических значений функции $f$. Подошва
$\d\DD_{[f]_\isot}$ ручки $\DD_{[f]_\isot}$ содержится в объединении
ручек $\DD_{[g]_\isot}$, таких что $[f]_\isot\prec[g]_\isot$.

{\rm (B)} Дискретная группа $\D/\D^0$ действует
на $\KK$ автоморфизмами косого цилиндрически-поли\-эд\-раль\-ного
комплекса, причем индуцированное действие на множестве ручек
согласовано с естественным действием на множестве $F^1/\sim_\isot$
классов изотопности функций.
В частности, для любого класса изотопности $[f]_\isot$ все ручки
$\DD_{[fh]_\isot}$, $h\in\D$, изоморфны одной и той же стандартной
косой цилиндрической ручке $(D_{[f]}\times\SS_{[f]})/\Gamma_{[f]}$,
см.\ определение~{\rm\ref {def:thick:cylinder}(В)}. Имеется
$\D/\D^0$-эквивариантный гомеоморфизм полиэдра $\KK$ на 
$\D/\D^0$-инвариантное подмножество некоторого гладкого $3q$-мерного
многообразия $\MM$ с плоской аффинной связностью, на котором группа
$\D/\D^0$ действует диффеоморфизмами, сохраняющими связность.

{\rm (C)} Для каждой ручки
$\DD_{[f]_\isot}\approx(D_{[f]}\times\SS_{[f]})/\Gamma_{[f]}\approx(D_{[f]}\times(\RR^{c([f])}\times(S^1)^{d([f])})\times
P_{[f]})/\Gamma_{[f]}\sim(S^1)^d/\Gamma_{[f]}$ размерность $d=d([f])$
тора $(S^1)^d$ обладает свойствами {\rm(\ref {eq:3.23.5})},
$c+d=n([f])$ {\rm(см.\ (\ref {eq:Zj}))} и
$d\le\min\{p'+p''+r'+r'',t-1\}$, где $t=t([f])\le q$ --- количество
связных компонент графа $G_f$ {\rm(см.\ обозначение~\ref
{not:Gf:or})}. Если число фиксированных критических точек
$p^*+q^*+r^*\le\chi(M)+1$, то $d=t-1$, а при дополнительном условии
$t=q$ выполнено $d=p'+p''+r'+r''$.
 \end{Thm}

 \begin{Rem} \label{rem:thm}
Согласно~\cite {kp1}, пространства $F^1,\FF^1$ суть сильные
деформационные ретракты пространств $F,\FF$ соответственно, а
забывающие отображения $\FF\to F$, $\FF^1\to F^1$ суть гомотопические
эквивалентности, где $\FF^1\subset\FF$ -- пространства оснащенных
функций Морса. Можно показать, что многообразие $\MM$ из теоремы~\ref
{thm:KP4}(B) в действительности гомеоморфно $\FF^1/\D^0$ (т.е.\
является универсальным пространством модулей оснащенных функций
Морса), причем действие группы $\D^0$ на $\FF^1$ свободно, проекция
$\FF^1\to\FF^1/\D^0$ является тривиальным расслоением со слоем
$\D^0$, а $\KK$ есть сильный деформационный ретракт $\MM$. Отсюда и
из (\ref {eq:EE}) следует требуемая гомотопическая эквивалентность
(\ref {eq:sim}). Можно также показать, что ограничения указанных
гомотопических эквивалентностей на любой класс $[f]_\isot$
изотопности функций из $F^1$ или на любую косую ручку комплекса $\KK$
являются гомотопическими эквивалентностями.
 \end{Rem}

Пусть $\Bbbk$ -- поле (например, $\RR,\QQ$ или $\ZZ_p$). Для
топологического пространства $X$ рассмотрим числа Бетти
$\beta_j(X):=\dim_\Bbbk H_j(X;\Bbbk)$ и полином Пуанкаре
$P(X,t):=\sum_{j=0}^\infty t^j\beta_j(X)$. Следующее утверждение
выводится из теоремы \ref {thm:KP4} стандартными методами теории
Морса (см., например,~\cite[\S45]{FF}).

\begin{Cor} \label {cor:ineq}
{\rm(A)} Если
количество $\hat p+\hat q+\hat r$ отмеченных критических точек
превосходит $\chi(M)$, то $\beta_j(\KK)=0$ при любом $j\ge3q-2$.

{\rm(B)} Пусть $\bar M=S^2$ {\rm(см.\ обозначение \ref {not:R:G0})},
$p^*+q^*+r^*\le\chi(M)+1\le\hat p+\hat q+\hat r$. Тогда $\D=\D^0$,
комплекс $\KK=\IKK$ является конечным, связным и компактным косым
торически-полиэдральным комплексом ранга $q-1$ и размерности $3q-2$
или $0$ (при $q\ge2$ и $q\le1$ соответственно);
числа Бетти $\beta_j=\beta_j(\IKK)$ комплекса $\IKK$ удовлетворяют
неравенствам Морса-Смейла:
 $$
 \beta_j-\beta_{j-1}+\beta_{j-2}-\beta_{j-3}+\ldots\le
 q_j-q_{j-1}+q_{j-2}-q_{j-3}+\ldots, \qquad j\ge0,
 $$
где $Q(t)=\sum_{j=0}^\infty t^jq_j:=\sum_{[f]\in
F^1/\sim}t^{q-s(f)}P(\DD_{[f]},t)$.
В частности, справедливы неравенства Морса:
 $$
 \chi(\IKK)=(-1)^{q-1}\left|\left\{[f]\in F^1/\sim \ \mid\
s(f)=1\right\}\right|, \qquad \beta_j\le q_j, \quad j\ge0.
 $$
\end{Cor}

\begin{Rem} \label {rem:numeration}
Всюду в формулировках утверждений настоящей статьи рассматриваются
произвольные обобщенные пространства функций Морса
$F=F_{p^*,p',p'';q^*,q',q'';r^*,r',r''}(M,\d^+M,\d^-M)$ и $F^1$,
состоящие из функций Морса,
у которых могут быть как отмеченные (пронумерованные), так и
неотмеченные (непронумерованные) критические точки. Однако в
доказательствах иногда будем считать, что $F=F^\num$ и
$F^1=F^{1,\num}$
(см.\ определение~$\ref {def:F}$), т.е.\ что все критические точки
функций Морса $f\in F$ пронумерованы. Это не ограничивает общности,
так как все отображения, построенные в настоящей статье, согласованы
с перенумерациями тех критических точек, которые изначально не были
отмечены (т.е.\
$\Sigma_{p''}\times\Sigma_{q''}\times\Sigma_{r''}$-эквивариантны
относительно действия группы
$\Sigma_{p''}\times\Sigma_{q''}\times\Sigma_{r''}$ на пространстве
$F^\num$ перенумерациями изначально неотмеченных критических точек).
Аналогично определению \ref {def:F} и обозначению \ref {not:R:G0}
обозначим через
$$
 \N_f := \N_{f,0} \cup \N_{f,1} \cup \N_{f,2},
 \qquad \hat\N_f:=\hat\N_{f,0}\cup \hat\N_{f,1}\cup \hat\N_{f,2}
 $$
множество всех критических точек (соответственно всех отмеченных
критических точек) функции $f\in F$. Имеем включения
$\N\subseteq\hat\N_{f}\subseteq\N_{f}$ и
$\N_\lam\subseteq\hat\N_{f,\lam}\subseteq\N_{f,\lam}$ множеств
фиксированных критических точек, отмеченных критических точек и всех
критических точек (соответственно индекса $\lambda$) функции $f$,
$\lam=0,1,2$.
\end{Rem}

\section{Построение стандартных косых цилиндрических ручек $\DD_{[f]_\isot}^\st$ и отображений
инцидентности $\chi_{[f]_\isot,[g]_\isot}$} \label {sec:KK0}

В данном параграфе предполагается, что выполнено условие (\ref
{eq:main}) (т.е.\ количество $\hat p+\hat q+\hat r$ отмеченных
критических точек превосходит $\chi(M)$). Для каждого класса
изотопности $[f]_\isot$ функций Морса мы опишем построение
стандартной косой цилиндрической ручки $\DD_{[f]_\isot}^\st$, а для
каждой пары примыкающих классов изотопности -- соответствующее
отображение инцидентности. В \S\ref {sec:KK} будет описано построение
``комплекса оснащенных функций Морса'' $\KK$, полученного из
описанных ручек при помощи описанных отображений инцидентности.
Проведем построение
в несколько шагов.

\subsection{Построение многогранника $D_{[f]_\isot}$ для класса изотопности $[f]_\isot$}\label{subsec:polyheder}

{\it Шаг 1} (определение пермутоэдра $\P^{q-1}$ порядка $q$ и
описание его граней). {\it Пермутоэдр порядка $q$} --- это выпуклый ($q-1$)-мерный
многогранник $\P^{q-1}$, вложенный в $q$-мерное пространство, вершины
которого получены перестановками координат вектора $(1,\dots,q)$
(впервые такие многогранники изучал Schoute (1911), название
появилось в книге Guiband \& Rosenstiehl (1963), более общие
``перестановочные многогранники'' с множеством вершин $\Sigma_q$
изучены Bowman (1972), см.\ также~\cite[доказательство теоремы 3, шаг
1]{BaK}). Опишем его подробнее: пусть $e_1,\dots,e_q$ -- стандартный базис $\RR^q$, и пусть $\P^{q-1}\subset\RR^q$ --
выпуклая оболочка множества точек $P_\pi=\sum_{k=1}^q
\left(k-\frac{q+1}{2}\right)e_{\pi_k}$, $\pi\in\Sigma_q$. Известно~\cite{P}, что
$\P^{q-1}$ -- $(q-1)$-мерный выпуклый многогранник в евклидовом
пространстве $\EE^{q-1}:=(e_1+\ldots+e_q)^\perp$, имеющий ровно $q!$
вершин $P_\pi$, $\pi\in\Sigma_q$, причем его $(q-s)$-мерные грани
находятся во взаимно однозначном соответствии с упорядоченными
разбиениями $J=(J_1,\dots,J_s)$ множества $\{1,\dots,q\}$ на $s$
непустых подмножеств $J_1,\dots,J_s$ (т.е.\ $\{1,\dots,q\}=J_1\sqcup
\ldots \sqcup J_s$), $1\le s\le q$. А именно, грань
$\tau_J\subset\P^{q-1}$, отвечающая разбиению $J=(J_1,\dots,J_s)$, --
это выпуклая оболочка множества точек
$(\Sigma_{r_1}\times\Sigma_{r_2-r_1}\times\ldots\times\Sigma_{r_s-r_{s-1}})(P_\pi)$,
где числа $0=r_0<r_1<\ldots<r_{s-1}<r_s=q$ и перестановка
$\pi\in\Sigma_q$ однозначно определяются условиями
 \begin{equation}\label{eq:J}
J_1=\{\pi_1,\dots,\pi_{r_1}\}, \ J_2=\{\pi_{r_1+1},\dots,\pi_{r_2}\},
\ \ldots, \ J_s=\{\pi_{r_{s-1}+1},\dots,\pi_{r_s}\},
 \end{equation}
$\pi_1<\ldots<\pi_{r_1}$, $\pi_{r_1+1}<\ldots<\pi_{r_2}, \ldots,
\pi_{r_{s-1}+1}<\ldots<\pi_{r_s}$. Здесь
$\Sigma_{r_1}\times\Sigma_{r_2-r_1}\times\ldots\times\Sigma_{r_s-r_{s-1}}$
-- подгруппа группы $\Sigma_q$, отвечающая разбиению
$\{1,\dots,q\}=\{1,\dots,r_1\}\sqcup\{r_1+1,\dots,r_2\}\sqcup\ldots\sqcup\{r_{s-1}+1,\dots,r_s\}$,
и действие перестановки $\rho\in\Sigma_q$ на точке $P_\pi$ дает точку
$P_{\rho\pi}$, где $(\rho\pi)_i:=\pi_{\rho_i}$, $1\le i\le q$.

Упорядоченные разбиения $J=(J_1,\dots,J_s)$ множества $\{1,\dots,q\}$
можно рассматривать как отношения частичного порядка на множестве
$\{1,\dots,q\}$. Если разбиение $\hat J$ получается из разбиения
$J=(J_1,\dots,J_s)$ путем измельчения (т.е.\ разбиения некоторых
множеств $J_k$ на несколько подмножеств), будем писать $\hat J\prec
J$. Из описания граней многогранника $\P^{q-1}$ следует, что условие
$\hat J\prec J$ равносильно $\tau_{\hat J}\prec\tau_J$, т.е.\
примыканию граней (см.\ определение~\ref {def:pol}(D)).

\begin{Lem} [о гранях пермутоэдра $\P^{q-1}$] \label{lem:P}
Пусть фиксирована грань $\hat\tau\prec\P^{q-1}$. Для любой грани
$\tau\prec\P^{q-1}$, такой что $\hat\tau\prec\tau$, рассмотрим
соответствующее разбиение $J=(J_1,\dots,J_s)$ и последовательность
чисел $(|J_1|,\dots,|J_s|)$. Тогда сопоставление
$\tau\mapsto(|J_1|,\dots,|J_s|)$ (для
$\hat\tau\prec\tau\prec\P^{q-1}$) инъективно. В частности, любой
автоморфизм пермутоэдра $\P^{q-1}\subset\RR^q$, индуцированный
перестановкой координатных осей, допустим {\rm(см.\ определение~\ref
{def:thick:cylinder}(A))}.
\end{Lem}

\begin{proof}{} Пусть $\hat\tau=\tau_{\hat J}$, $\hat J=(\hat J_1,\dots,\hat
J_{\hat s})$. Ввиду $\hat J\prec J$ упорядоченное разбиение $J$
получается из упорядоченного разбиения $\hat J$ путем объединения
некоторых соседних подмножеств в одно подмножество, т.е.\
$J_k=\bigcup\limits_{i=a_{k-1}+1}^{a_k}\hat J_i$, $1\le k\le s$, для
некоторых $a_0=0<a_1<\ldots<a_s=\hat s$. Поэтому
$|J_k|=\sum\limits_{i=a_{k-1}+1}^{a_k}|\hat J_i|$. Отсюда следует,
что по разбиению $\hat J$ и набору чисел $(|J_1|,\dots,|J_s|)$
последовательность $a_0=0<a_1<\ldots<a_s=\hat s$, а потому и
разбиение $J$, определяются однозначно. Лемма доказана.
\end{proof}

{\it Шаг 2.} Для каждой функции Морса $f\in F$ рассмотрим множество
$\N_{f,1}=:\{y_j\}_{j=1}^q\approx\{1,\dots,q\}$ ее седловых
критических точек (см.\ замечание~\ref {rem:numeration}) и евклидово
векторное пространство $0$-коцепей
 \begin{equation} \label{eq:H0f}
 H_f^0:=C^0(\N_{f,1};\RR)=\RR^{\N_{f,1}}\cong\RR^q
 \end{equation}
со стандартной евклидовой метрикой. В этом векторном пространстве
рассмотрим многогранник $\P^{q-1}_f\subset H_f^0$, являющийся образом
многогранника $\P^{q-1}\subset\RR^q$ при какой-либо биекции
$\N_{f,1}\to\{1,\dots,q\}$. Рассмотрим ``вычисляющую'' 0-коцепь
 $$
\barc=\barc(f):=f|_{\N_{f,1}}=(c_1,\dots,c_q)\in H_f^0,
 $$
т.е.\ функцию $\barc\:\N_{f,1}\to\RR$, сопоставляющую каждой седловой
точке $y_j\in \N_{f,1}$ значение $c_j:=f(y_j)$ функции $f$ в этой
точке, $1\le j\le q$. Сопоставим 0-коцепи $\barc=(c_1,\dots,c_q)$
число $s(\barc):=|\{c_1,\dots,c_q\}|$ различных седловых значений и
упорядоченное разбиение $J=J(\barc)=(J_1,\dots,J_s)$ множества
седловых точек $\N_{f,1}\approx\{1,\dots,q\}$, определяемое
свойствами~(\ref {eq:J}) и $c_{\pi_1}=\ldots=c_{\pi_{r_1}} <
 c_{\pi_{r_1+1}}=\ldots=c_{\pi_{r_2}} < \ldots <
c_{\pi_{r_{s-1}+1}}=\ldots=c_{\pi_{r_s}}$. (То есть, $J$ -- это
отношение частичного порядка на множестве $\N_{f,1}$ седловых
критических точек функции $f$ значениями функции $f|_{\N_{f,1}}$.)
Можно также рассматривать $J=J(\barc)$ как сюръекцию
$\N_{f,1}\to\{1,\dots,s\}$, переводящую $y_{\pi_j}\mapsto k$ при
$r_{k-1}<j\le r_k$, $1\le j\le q$.

В каждом классе эквивалентности $[f]\in F^1/\sim$ (соответственно
классе изотопности $[f]_\isot\in F^1/\sim_\isot$) отметим ровно одну
функцию Морса $f$ этого класса, так чтобы любая функция $f\in F^1$,
являющаяся отмеченной функцией класса эквивалентности $[f]$, являлась
отмеченной функцией класса изотопности $[f]_\isot$. Сопоставим классу
изотопности $[f]_\isot$ с отмеченной функцией $f$ и разбиению
$J(\barc(f))$ грань
 $$
 D_{[f]_\isot} = D_f := \tau_{J(\barc(f))}\subset\P^{q-1}_f.
 $$

{\it Шаг 3.} Изучим взаимосвязь многогранников
$D_{[f]_\isot},D_{[g]_\isot}$ для примыкающих классов изотопности
$[f]_\isot\prec[g]_\isot$. Для любой функции $f\in F$ и
соответствующего $q$-мерного евклидова пространства $H^0_f\cong\RR^q$
(см.\ шаг 2) выполнены следующие два свойства:

1) для любого вектора $\barc\in H^0_f\cong\RR^q$ существует
$\eps_0>0$, такое что (i) для любого $\barc'\in H^0_f\cong\RR^q$ со
свойством $|\barc'-\barc|<\eps_0$ выполнено $J(\barc')\preceq
J(\barc)$, и (ii) для любых $\eps\in(0,\eps_0]$ и разбиения $\hat
J\preceq J(\barc)$ существует $\barc'\in\RR^q$ со свойствами
$|\barc'-\barc|<\eps$ и $J(\barc')=\hat J$;

2) согласно~\cite[утверждение~1.1 и~\S3]{K}, любая
(``невозмущенная'') функция $f\in F$ имеет столь малую окрестность
$\UU_f$ в $F$, что для любых (``возмущенных'') функций $\tilde
f,\tilde f_1\in\UU_f$ равенства $[\tilde fh_{0;f,\tilde
f}^{-1}]_\isot^\fix=[\tilde f_1h_{0;f,\tilde f_1}^{-1}]_\isot^\fix$ и
$J((h_{0;f,\tilde f}^{-1})^{*0}(\barc(\tilde f)))=J((h_{0;f,\tilde
f_1}^{-1})^{*0}(\barc(\tilde f_1)))$ равносильны, где через
$h_{0;f,\tilde f}\in\D^0$ обозначен диффеоморфизм, близкий к
тождественному, такой что $h_{0;f,\tilde f}^{-1}(\N_f)=\N_{\tilde
f}$, через $h_{0;f,\tilde f}^{*0}\:H_f^0\to H_{\tilde f}^0$
индуцированный изоморфизм групп $0$-коцепей, а через $[\tilde
fh_{0;f,\tilde f}^{-1}]_\isot^\fix$ обозначен класс изотопности
функции $\tilde fh_{0;f,\tilde f}^{-1}$ с фиксированным множеством
критических точек (при фиксированной функции $f$); в частности, при
выполнении указанных равносильных равенств существует диффеоморфизм
$h_{1;\tilde fh_{0;f,\tilde f}^{-1},\tilde f_1h_{0;f,\tilde
f_1}^{-1}}\in\D^0$,
гомотопный $\id_M$ в пространстве гомеоморфизмов пары $(M,\N_f)$ и
переводящий линии уровня функции $\tilde f_1h_{0;f,\tilde f_1}^{-1}$
в линии уровня функции $\tilde fh_{0;f,\tilde f}^{-1}$ с сохранением
направления роста.

В силу этих свойств, для любой функции $f\in F^1$ имеется сюръекция
$\delta[f]_\isot$ множества всех граней
$\tau'\prec\tau:=\tau_{J(\barc(f))}$ на множество всех классов
изотопности $[g]_\isot\succ[f]_\isot$ (см.\ определение~\ref
{def:pol}(D)), такая что
$\delta{[f]_\isot}\:\tau'\mapsto\delta_{\tau'}[f]_\isot:=[\tilde
f]_\isot$ тогда и только тогда, когда
 \begin{equation} \label {eq:incid}
 \tau'
 =\tau_{J((h_{0;f,\tilde f}^{-1})^{*0}(\barc(\tilde f)))}
 =\tau_{J((h_{f,\tau'}^{-1})^{*0}(\barc(g)))},
 \end{equation}
где $\tilde f\in\UU_f$ и $h_{0;f,\tilde f}$ как во втором свойстве
выше, $g$ --- отмеченная функция класса изотопности $[\tilde
f]_\isot$,
 \begin{equation} \label {eq:incid'}
 h_{f,\tau'}:=h_{0;f,\tilde f}h_{1;\tilde f,g}\in\D^0, \qquad
 h_{f,\tau'}^{*0}\: H_f^0\to H_g^0
 \end{equation}
--- индуцированный изоморфизм, а диффеоморфизм $h_{1;\tilde
f,g}\in\D^0$ переводит линии уровня функции $g$ в линии уровня
функции $\tilde f$ с сохранением направления роста (он существует
ввиду изотопности функций $\tilde f,g$), откуда
 $
\tilde f=h_2gh_{1;\tilde f,g}^{-1}
 $
для некоторого $h_2\in\Diff^+[-1;1]$. Корректность определения
сюръекции $\delta[f]_\isot$, т.е.\ независимость класса изотопности
$\delta_{\tau'}[f]_\isot=[\tilde f]_\isot$ от выбора функции $\tilde
f\in\UU_f$ с заданным значением $J((h_{0;f,\tilde
f}^{-1})^{*0}(\barc(\tilde f)))$, следует из второго свойства (см.\
выше). Если диффеоморфизм $\tilde h_{f,\tau'}=h_{0;f,\tilde
f_1}h_{1;\tilde f_1,g}$ построен с помощью функции $\tilde
f_1\in\UU_f$, такой что $J((h_{0;f,\tilde f}^{-1})^{*0}(\barc(\tilde
f)))=J((h_{0;f,\tilde f_1}^{-1})^{*0}(\barc(\tilde f_1)))$, то
диффеоморфизм $h_{f,\tau'}^{-1}h_{1;\tilde fh_{0;f,\tilde
f}^{-1},\tilde f_1h_{0;f,\tilde f_1}^{-1}}\tilde h_{f,\tau'}$ (см.\
второе свойство выше) переводит линии уровня функции $g$ в линии
уровня функции $g$ с сохранением направления роста, а потому ввиду
\cite [лемма 1]{KP2} принадлежит $(\stab_{\D^0}g) \
(\Diff^0(M,\N_g))$, откуда
 \begin{equation} \label {eq:corr}
 h_{f,\tau'}^{-1} \tilde h_{f,\tau'}
\in (\stab_{\D^0}g) \ (\Diff^0(M,\N_g)) ,
 \end{equation}
где через $\stab_{\D^0}g$ обозначена группа изотропии элемента $g$
относительно естественного правого действия группы $\D^0$ на $F^1$, а
через $\Diff^0(M,\N_{g})\subset\D^0$ группа диффеоморфизмов пары
$(M,\N_{g})$, гомотопных $\id_M$ в классе гомеоморфизмов пары, где
$\N_{g}$ --- множество критических точек функции $g$.

Пусть $f\in F^1$ --- отмеченная функция Морса класса изотопности
$[f]_\isot$. Для любой грани
$\tau'\prec\tau_{J(\barc(f))}=:D_{[f]_\isot}$ обозначим через $g\in
F^1$ отмеченную функцию класса изотопности $\delta_{\tau'}[f]_\isot$
и фиксируем диффеоморфизм $h_{f,\tau'}\in\D^0$ как в~(\ref
{eq:incid}) и~(\ref {eq:incid'}). Тогда, ввиду равенств~(\ref
{eq:incid}) и
$(h_{f,\tau'}^{-1})^{*0}(\tau_{J(\barc)})=\tau_{J((h_{f,\tau'}^{-1})^{*0}(\barc))}$,
имеем изоморфизм граней
 \begin{equation}\label {eq:glue:D}
 h_{f,\tau'}^{*0}|_{\tau'}\: \tau' \longrightarrow h_{f,\tau'}^{*0}(\tau')
 = \tau_{J(\barc(g))} = D_{[g]_\isot}.
 \end{equation}

Из указанных в начале шага двух свойств мы также получаем, что из
$[h]_\isot\succ[g]_\isot\succ[f]_\isot$ следует
$[h]_\isot\succ[f]_\isot$, и
 \begin{equation} \label {eq:*delta}
 \delta_{\tau''}{[f]_\isot}=\delta_{h_{f,\tau'}^{*0}(\tau'')}\delta_{\tau'}{[f]_\isot}
 \quad \mbox{для любых граней }\tau''\prec\tau'\prec\tau_{J(\barc(f))}.
 \end{equation}
Пусть $g,g_1\in F^1$ --- отмеченные функции классов изотопности
$\delta_{\tau'}[f]_\isot$, $\delta_{\tau''}[f]_\isot$ соответственно,
и пусть $\tilde g\in\UU_g$ и
$h_{f,\tau'}^{*0}(\tau'')=\tau_{J((h_{0;g,\tilde
g}^{-1})^{*0}(\barc(\tilde g)))}$ (см.~(\ref {eq:incid}), (\ref
{eq:*delta})). Покажем, что выполнен следующий аналог {\it
соотношения транзитивности}:
 \begin{equation} \label {eq:obstruction}
 h_{f,\tau''}^{-1} h_{f,\tau'}h_{g,h_{f,\tau'}^{*0}(\tau'')} \in  (\stab_{\D^0}g_1) \ (\Diff^0(M,\N_{g_1})) ,
 \end{equation}
где $\N_{g_1}$ --- множество критических точек функции $g_1$.
Действительно, функция $\tilde f_1:=h_2\tilde gh_{1;\tilde f,g}^{-1}
\in [g_1]_\isot$ близка к функции $\tilde f=h_2gh_{1;\tilde
f,g}^{-1}$ (а потому и к функции $f$); диффеоморфизм $h_{1;\tilde
f,g}h_{0;g,\tilde g}h_{1;\tilde f,g}^{-1}$ близок к $\id_M$ и
переводит критические точки ``возмущенной'' функции $\tilde f_1$ в
критические точки ``невозмущенной'' функции $\tilde f$, а потому
диффеоморфизм
 $$
  \tilde h_{0;f,\tilde f_1}:= h_{0;f,\tilde f} \ h_{1;\tilde f,g}h_{0;g,\tilde g}h_{1;\tilde f,g}^{-1} \in \D^0
 $$
близок к $\id_M$ и переводит критические точки ``возмущенной''
функции $\tilde f_1$ в критические точки ``невозмущенной'' функции
$f$; диффеоморфизм $h_{1;\tilde f,g}\in\D^0$ переводит линии уровня
функции $\tilde g$ в линии уровня функции $\tilde f_1$ с сохранением
направления роста, а потому диффеоморфизм
 $$
 \tilde h_{1;\tilde f_1,g_1}:=h_{1;\tilde f,g}h_{1;\tilde g,g_1} \in\D^0
 $$
переводит линии уровня функции $g_1$ в линии уровня функции $\tilde
f_1$ с сохранением направления роста. Отсюда
 $$
   h_{f,\tau'}h_{g,h_{f,\tau'}^{*0}(\tau'')}
 = h_{0;f,\tilde f}h_{1;\tilde f,g}\ h_{0;g,\tilde g}h_{1;\tilde g,g_1}
 =  \tilde h_{0;f,\tilde f_1}                                                \tilde h_{1;\tilde f_1,g_1},
 $$
т.е.\ мы разложили диффеоморфизм
$h_{f,\tau'}h_{g,h_{f,\tau'}^{*0}(\tau'')}$ в композицию, аналогичную
разложению $h_{f,\tau''}=h_{0;f,\tilde f_1}h_{1;\tilde f_1,g_1}$,
см.~(\ref {eq:incid'}). Ввиду (\ref {eq:corr}) это доказывает (\ref
{eq:obstruction}).

Если $h\in\stab_{\D^0}f$, а $\tau',g$ как выше, то ввиду (\ref
{eq:incid}) выполнено $h^{*0}(\tau')=\tau_{J((h_{0;f,\tilde
fh}^{-1})^{*0}(\barc(\tilde fh)))}$ и
$\delta_{h^{*0}(\tau')}[f]_\isot=[gh]_\isot=[g]_\isot$ (ввиду
$h\in\D^0$), откуда $h_{f,h^{*0}(\tau')}=h_{0;f,\tilde fh}h_{1;\tilde
fh,g} 
\\
=h_0h^{-1}h_{0;f,\tilde f}h\ h^{-1}h_{1;\tilde f,g}
 h_1
=h_0h^{-1}h_{f,\tau'} h_1$ для некоторых $h_0\in\Diff^0(M,\N_f)$ и
$h_1\in\D^0$, таких что $h_1$ переводит линии уровня функции $g$ в
линии уровня функции $g$ с сохранением направления роста, поэтому
$h_1\in(\stab_{\D^0}g) (\Diff^0(M,\N_g))$ ввиду \cite [лемма 1]{KP2},
откуда
 \begin{equation} \label {eq:stab}
 h_{f,\tau'}^{-1} h h_{f,h^{*0}(\tau')}
\in (\stab_{\D^0}g) \ (\Diff^0(M,\N_g)) .
 \end{equation}

\subsection{Построение утолщенного цилиндра $\SS_{[f]_\isot}$ для класса изотопности $[f]_\isot$}\label{subsec:cylinder}

{\it Шаг 4.} Пусть $f\in F^1$ --- отмеченная функция Морса класса
изотопности $[f]_\isot$. По аналогии с пространством 0-коцепей
$H^0_f\cong\RR^q$ (см.~(\ref {eq:H0f})) введем двойственные друг
другу векторные пространства относительных 1-гомологий и
относительных 1-когомологий над полем $\RR$:
 \begin{equation} \label {eq:H1f}
 \begin{array}{l}
 H_{f,1}:=H_1(M\setminus(\N_{f,0}\cup \N_{f,2}),\N_{f,1};\RR)\cong\RR^{2q}, \\
 \phantom{0} H^1_f:=H^1(M\setminus(\N_{f,0}\cup \N_{f,2}),\N_{f,1};\RR) \cong \Hom_\RR(H_{f,1},\RR)
 \cong \RR^{2q},
 \end{array}
 \end{equation}
где изоморфизм $H^1_f\cong\Hom_\RR(H_{f,1},\RR)$ индуцирован
равенством $C^q(X,A;\RR)=\Hom(C_q(X)/C_q(A),\RR)$ для пары Борсука
$(X,A):=(M\setminus(\N_{f,0}\cup \N_{f,2}),\N_{f,1})$ при $q=1$
(см.~\cite[\S\S12.5, 15.2, 15.5]{FF}). Рассмотрим ориентированный
граф $G_f\subset M\setminus(\N_{f,0}\cup \N_{f,2})$, см.\
обозначение~\ref {not:Gf:or}. Этот граф имеет $q$ вершин (являющихся
седловыми точками $\b_1,\dots,\b_q\in \N_{f,1}$), степени всех вершин
равны 4, а значит в графе $2q$ ребер, которые обозначим
$e_1,\dots,e_{2q}$. Обозначим относительный гомологический класс
ориентированного ребра $e_i$ через $[e_i]\in H_{f,1}$, $1\le i\le2q$.

Определим в векторном пространстве $H^1_f\cong\RR^{2q}$ выпуклые
подмножества
 \begin{equation} \label {eq:Us}
 U_{[f]_\isot}
 \subset U_{[f]_\isot}^{\infty} \subset H^1_f
 \end{equation}
системами из $4q$ и $2q$ неравенств соответственно:
 \begin{equation} \label {eq:U}
 U_{[f]_\isot}=U_f:=\left\{u\in H^1_f \ \left|\ 1\le u([e_i]) \le \frac {(2q-1)!!}{(2q-2s+1)!!},\ 1\le i\le 2q \right. \right\},
 \end{equation}
 \begin{equation} \label {eq:U0}
 U_{[f]_\isot}^{\infty}=U_f^{\infty}:=\left\{u\in H^1_f \ \left|\ u([e_i]) >0,\ 1\le i\le 2q \right. \right\},
 \end{equation}
где $s=s(\barc(f)):=|f(\N_{f,1})|$ -- количество седловых значений функции $f$.

{\it Шаг 5.} Каждая компонента связности пространства $M\setminus
G_f$ содержит не более одной критической точки функции $f$ (а именно,
точки минимума или максимума) и не более одной компоненты связности
края $\d M$, и гомеоморфна либо открытому кругу (с одной критической
точкой), либо полуоткрытому цилиндру $S^1\times[-1;0)$ или
$S^1\times(0;1]$ (c одной компонентой края
$S^1\times\{\pm1\}\subset\d^\pm M$), либо ``открытому цилиндру''
 \begin{equation}\label {eq:Zj}
Z_\ell=Z_\ell(f)\approx S^1\times (0;1), \quad 1\le\ell\le n=n(f),
 \end{equation}
которые вместе со своим замыканием содержатся в $(\Int
M)\setminus(\N_{f,0}\cup \N_{f,2})$, где $n=n(f)$ --- количество
открытых цилиндров. Сопоставим открытому цилиндру $Z_\ell$ его
серединную окружность
 \begin{equation} \label {eq:gamma:j}
 \gamma_\ell=\gamma_\ell(f)=S^1\times \left\{\frac12\right\} \subset Z_\ell=S^1\times(0;1)
 \end{equation}
и следующее линейное векторное поле $v_\ell$ на векторном
пространстве $H^1_f$. Значение $v_\ell(u)\in H^1_f$ поля $v_\ell$ в
любой точке $u\in H^1_f$ --- это относительный 1-коцикл, значение
которого на любом относительном 1-цикле $a\in
H_1(M\setminus(\N_{f,0}\cup \N_{f,2}),\N_{f,1})\subset H_{f,1}$
определяется формулой
 \begin{equation} \label {eq:v:ell}
 v_\ell(u)([a]) := \<[\gamma_\ell], a\> \ u([\gamma_\ell]), \quad u\in H^1_f,\ a\in
 H_1\left(M\setminus(\N_{f,0}\cup \N_{f,2}),\N_{f,1}\right),
 \end{equation}
где $\<[\gamma_\ell], a\>$ --- индекс пересечения цикла
$[\gamma_\ell]\in H_1(M\setminus \N_f) \subset H_{f,1}$ и
относительного цикла $a$, $1\le \ell\le n$. Другими словами, линейное
векторное поле $v_\ell$ на $H^1_f $ задается $\RR$-линейным
оператором $H^1_f\to H^1_f$, являющимся обратным образом при
изоморфизме $H^1_f\cong\Hom_\RR(H_{f,1},\RR)$ $\RR$-линейного
оператора $H_{f,1}\to H_{f,1}$,
$a\mapsto\<[\gamma_\ell],a\>[\gamma_\ell]$. На шаге 6 мы установим
свойства векторных полей $v_1,\dots,v_n$. Определим диффеоморфизм
 \begin{equation}\label {eq:hj}
h_\ell:=t_{\gamma_\ell} \in \stab_\D f \subset\D, \quad 1\le \ell\le
n,
 \end{equation}
как {\it скручивание Дэна} $t_{\gamma_\ell}$ (см.~\cite {Dehn})
вокруг окружности $\gamma_\ell$ (скручивание Дэна $t_{\gamma_\ell}$
совпадает с $\id_M$ вне открытого цилиндра $Z_\ell$ и получается с
помощью разрезания поверхности вдоль окружности $\gamma_\ell$,
перекручивания одного конца разреза на $2\pi$ и cклеивания).
Диффеоморфизмы $h_\ell$, $1\le\ell\le n$, попарно коммутируют, так
как их носители попарно не пересекаются. Рассмотрим порожденную ими
абелеву группу
$$
 \tildeGamma_{[f]_\isot}=\tildeGamma_f:=\<h_1,\dots,h_n\>\subset \stab_\D f\subset\D.
$$
Рассмотрим индуцированные автоморфизмы $h_\ell^*\in\Aut(H^1_f)$,
$1\le\ell\le n$, и порожденную ими абелеву группу
\begin{equation} \label {eq:Theta}
\tildeGamma_{[f]_\isot}^*=\tildeGamma_f=\<h_1^*,\dots,h_n^*\>
 \subset\Aut(H^1_f).
\end{equation}
Нетрудно показать, что подгруппа
$\tildeGamma_{f}^*\subset\Aut(H^1_f)$ изоморфна свободной абелевой
группе ранга $n$, и что автоморфизм $h_\ell^*$ совпадает с потоком
$g_{v_\ell}^1\:H^1_f\to H^1_f$ векторного поля $v_\ell$ за время 1,
$1\le\ell\le n$. Рассмотрим в группе $\tildeGamma_f^*\cong\ZZ^n$
подгруппу $(\D^0\cap\tildeGamma_f)^*\subset\tildeGamma_f^*$ и
рассмотрим пространства орбит
 $$
 \SS_{[f]_\isot} = \SS_f := U_f / (\D^0 \cap \tildeGamma_f)^*, \quad  \SS^\infty_{[f]_\isot} = \SS^\infty_f := U_f^\infty / (\D^0 \cap \tildeGamma_f)^*.
 $$
Рассуждения на следующих шагах проводятся для $U_f$, но верны и для
$U_f^\infty$.

{\it Шаг 6.} На этом шаге определяется свободное действие цилиндра
$\RR^{n-d}\times(S^1)^d$ на пространстве $\SS_f$, где $d=d([f])$ ---
ранг группы $(\D^0 \cap \tildeGamma_f)^*$ (как свободной абелевой
группы).

Построим явно базис векторного пространства $H^1_f\cong\RR^{2q}$.
Если количество $n=n(f)$ открытых цилиндров $Z_\ell\subset M\setminus
G_f$ (см.\ (\ref {eq:Zj})) равно нулю, положим $\tilde G_f:=G_f$.
Если $n>0$, будем выкидывать из графа $G_f$ по одному (открытому)
ребру, чтобы каждый раз количество компонент связности дополнения
графа в $M$, не пересекающихся с $(\d M)\cup \N_{f,0}\cup \N_{f,2}$,
уменьшалось на 1. Так как для графа $G_f$ указанное количество равно
$n$, то после выкидывания из него $n$ ребер (для определенности
$e_1,\dots,e_n$) указанным алгоритмом их количество станет равным
нулю, и получится подграф с $q$ вершинами и $2q-n$ ребрами
$e_{n+1},\dots,e_{2q}$. В каждом открытом цилиндре $Z_\ell$ проведем
(отрытое) ориентированное ребро $\tilde e_\ell$, гладко вложенное в
этот цилиндр, с концами в седловых точках, такое, что ограничение
функции $f$ на это ребро монотонно возрастает. Добавим к полученному
подграфу $n$ ориентированных ребер
 $$
\tilde e_\ell \subset Z_\ell, \qquad 1\le \ell\le n
 $$
(взамен выброшенных $e_1,\dots,e_n$). В результате получим граф
$\tilde G_f \subset (\Int M)\setminus(\N_{f,0}\cup \N_{f,2})$ с $2q$
ребрами $\tilde e_1,\dots,\tilde e_n$, $e_{n+1},\dots,e_{2q}$ и $q$
вершинами $\b_1,\dots,\b_q$. Так как дополнение графа $\tilde G_f$ в
поверхности $M$ состоит из открытых кругов (содержащих по одной точке
минимума или максимума) и полуоткрытых цилиндров (содержащих по одной
компоненте края $M$), то граф $\tilde G_f$ является сильным
деформационным ретрактом поверхности $M\setminus(\N_{f,0}\cup
\N_{f,2})$. Следовательно, вложение $\tilde G_f\hookrightarrow
M\setminus(\N_{f,0}\cup \N_{f,2})$ индуцирует изоморфизмы $2q$-мерных
векторных пространств
$$
 H_1(\tilde G_f,\N_{f,1};\RR)\cong H_{f,1}, \quad
 H^1(\tilde G_f,\N_{f,1};\RR)\cong H^1_f.
$$
Относительные классы гомологий $[\tilde e_1],\dots,[\tilde e_n]$,
$[e_{n+1}],\dots,[e_{2q}] \in H_1(\tilde G_f,\N_{f,1};\RR)$
ориентированных ребер графа $\tilde G_f$ образуют базис векторного
пространства $H_1(\tilde G_f,\N_{f,1};\RR)\cong
H_{f,1}\cong\RR^{2q}$. Переход к двойственному базису дает базис
$[\tilde e_1]^*,\dots,[\tilde e_n]^*$, $[e_{n+1}]^*,\dots,[e_{2q}]^*$
векторного пространства
\\
$H^1(\tilde G_f,\N_{f,1};\RR)\cong
H^1_f\cong(H_{f,1})^*\cong\RR^{2q}$.

Пусть $\tilde u_1,\dots,\tilde u_n,u'_{n+1},\dots u'_{2q}$ --
координаты в $H^1_f$ по отношению к базису $[\tilde
e_1]^*,\dots,[\tilde e_n]^*$, \\ $[e_{n+1}]^*,\dots,[e_{2q}]^*$.
Рассмотрим разложение
 \begin{equation}\label{eq:oplus}
 H^1_f=\<[\tilde e_1]^*,\dots,[\tilde e_n]^*\>\oplus\<[e_{n+1}]^*,\dots,[e_{2q}]^*\>.
 \end{equation}
Представим любой относительный коцикл $u\in H^1_f$ как сумму
$u=\tilde u+u'$ его проекций на подпространства в
разложении~(\ref{eq:oplus}). Нетрудно доказывается, что
\begin{equation} \label {eq:U'}
 [e_i]\in\<[e_{n+1}],\dots,[e_{2q}]\>=\Im\left[H_1(G_f,\N_{f,1};\RR)\to H_{f,1}\right], \quad 1\le i\le 2q,
 \end{equation}
$$
 \<[\tilde e_\ell]^*\>_{\ell=1}^n = \ker\left[H^1_f\to H^1(G_f,\N_{f,1};\RR)\right],
 \quad
 \<[e_{i}]^*\>_{i=n+1}^{2q}\cong
  H^1(G_f,\N_{f,1};\RR).
$$
Отсюда $u([e_i])=u'([e_i])$, $1\le i\le2q$. Положим
\begin{equation} \label {eq:U''}
U'_f:=\left\{u'\in\<[e_{i}]^*\>_{i=n+1}^{2q}\ \left|\ 1\le
u'([e_\ell])\le \frac {(2q-1)!!}{(2q-2s+1)!!},\ 1\le\ell\le
2q\right.\right\},
\end{equation}
ср.~(\ref {eq:U}). Тогда $u\in U_f$ в том и только том случае, когда
$u'\in U'_f$, причем $U'_f$ -- выпуклый многогранник. Поэтому
справедливо разложение
 \begin{equation}\label{eq:oplus'}
 U_f=\<[\tilde e_1]^*,\dots,[\tilde e_n]^*\>\oplus U'_f,
 \qquad \mbox{где}\quad U'_f\subset\<[e_{n+1}]^*,\dots,[e_{2q}]^*\>.
 \end{equation}

В базисе $[\tilde e_1]^*,\dots,[\tilde e_n]^*$,
$[e_{n+1}]^*,\dots,[e_{2q}]^*$ пространства $H^1_f$ линейные
векторные поля $v_1,\dots,v_n$ на $H^1_f$ имеют вид
\begin{equation} \label {eq:v}
 v_\ell(u)=u([\gamma_\ell])\ [\tilde e_\ell]^*, \quad
 \<v_\ell(u)\>_{\ell=1}^n\subseteq\<[\tilde e_\ell]^*\>_{\ell=1}^n
 = \ker\left[H^1_f\to H^1(G_f,\N_{f,1};\RR)\right],
\end{equation}
т.е.\ касательны каждой $n$-мерной плоскости $\<[\tilde
e_1]^*,\dots,[\tilde e_n]^*\>+u'\subset H^1_f$,
$u'\in\<[e_{n+1}]^*,\dots,[e_{2q}]^*\>$, и всюду на ней имеют
постоянные коэффициенты ввиду (\ref {eq:U'}). Поэтому каждая такая
$n$-мерная плоскость инвариантна относительно потоков векторных полей
$v_1,\dots,v_n$ на $H^1_f$ и эти поля коммутируют. Так как при $u'\in
U'_f$ эти векторные поля (с постоянными коэффициентами) линейно
независимы в указанной плоскости, то их потоки $g_{v_1}^{t_1}\ldots
g_{v_n}^{t_n}$ порождают свободное действие группы $\RR^n$ на
пространстве $U_f$, см.~(\ref {eq:oplus'}), причем орбиты этого
действия являются $n$-мерными плоскостями $\<[\tilde
e_1]^*,\dots,[\tilde e_n]^*\>+u'\subset U_f$, $u'\in U_f'$. Так как
группа $\RR^n$ действует свободно на $U_f$, то ее стандартная
целочисленная решетка $\ZZ^n\subset\RR^n$ тоже действует свободно на
$U_f$. Так как действие $\ell$-го базисного элемента решетки $\ZZ^n$
совпадает с $\ell$-ым базисным элементом
$g_{v_\ell}^1=h_\ell^*\in\Aut(H^1_f)$ группы
$\tildeGamma_f^*\cong\ZZ^n$ (см.\ (\ref {eq:hj}), (\ref {eq:Theta})),
то действие группы $\tildeGamma_f^*\subset\Aut(H^1_f)$ на $U_f$
свободно и коммутирует с действием $\RR^n$ на $U_f$ (заданным при
помощи потоков векторных полей $v_1,\dots,v_n$). Поэтому действие
группы $\RR^n$ на $U_f$ индуцирует корректно определенное свободное
действие цилиндра $\RR^n/Z^d\cong\RR^{n-d}\times(S^1)^d$ на
факторпространстве $\SS_f = U_f / (\D^0 \cap \tildeGamma_f)^*$, где
$(\D^0 \cap \tildeGamma_f)^*\cong Z^d\subset\ZZ^n$ --- подгруппа
группы $\tildeGamma_f^*\cong\ZZ^n\subset\RR^n$, и через $d=d([f])$
обозначен ее ранг (как ранг свободной абелевой группы).

Все рассуждения данного шага верны для $U_f^\infty,(U_f')^\infty$
вместо $U_f,U_f'$, где $(U_f')^\infty$ определяется аналогично (\ref
{eq:U''}).

{\it Шаг 7.} На этом шаге вводится на пространстве $\SS_f$ структура
утолщенного цилиндра (см.\ определение~\ref {def:thick:cylinder}).
Для этого будут построены специальные (криволинейные) координаты на
выпуклом множестве $U_f\subset H^1_f$ и на утолщенном цилиндре
$\SS_f=U_f/(\D^0 \cap \tildeGamma_f)^*$, в которых построенные выше
свободные действия группы $\RR^n$ на $U_f$ и цилиндра
$\RR^n/Z^d\cong\RR^{n-d}\times(S^1)^d$ на $\SS_f$ ``выпрямляются''.

Построим явно набор образующих группы
$(\D^0\cap\tildeGamma_f)^*\subset\tildeGamma_f^*\subset\Aut(H^1_f)$.
Напомним, что набором свободных образующих группы
$\tildeGamma_f^*\cong\ZZ^n$ является набор автоморфизмов
$h_1^*,\dots,h_n^*\in\Aut(H^1_f)$, отвечающих открытым цилиндрам
$Z_1,\dots,Z_n$, где $n=n(f)$, см.~(\ref {eq:Zj}). Покажем, что после
подходящей перенумерации цилиндров $Z_\ell$ (и отвечающих им
автоморфизмов $h_\ell^*$) подгруппа
$(\D^0\cap\tildeGamma_f)^*\subset\tildeGamma_f^*\cong\ZZ^n$
раскладывается в прямое произведение подгрупп
 \begin{equation}\label{eq:generH}
 (\D^0\cap\tildeGamma_f)^*
 = \tildeGamma^*_{f,0} \times \tildeGamma^*_{f,1} \times
\ldots \times \tildeGamma^*_{f,e}, \qquad \mbox{где}\quad
 \tildeGamma^*_{f,0}=\<h_1^*,\dots,h_{\nu_0}^*\>,
\end{equation}
 $$ \tildeGamma^*_{f,k}=\<
 (h^*_{\nu_{k-1}+1})^{-1}h^*_{\nu_{k-1}+2},\
 (h^*_{\nu_{k-1}+2})^{-1}h^*_{\nu_{k-1}+3},\ \dots\ ,\
 (h^*_{\nu_k-1})^{-1}    h^*_{\nu_k}\>, \
1\le k\le e,
 $$
где целые числа $0\le e\le n$ и
$0=\nu_{-1}\le\nu_0<\nu_1<\ldots<\nu_{e}\le n$ зависят от $[f]$.
Отсюда следует, что ранг группы $(\D^0\cap\tildeGamma_f)^*$ равен
 \begin{equation} \label {eq:3.23.5}
 \rank(\D^0\cap\tildeGamma_f)^*=d=\nu_e-e.
 \end{equation}
Из (\ref {eq:3.23.5}) нетрудно вывести, что он не превосходит числа
$p'+p''+r'+r''$ ``плавающих'' точек локальных минимумов и максимумов,
а также получить остальные оценки для $d$ из теоремы~\ref
{thm:KP4}(C).

{\it Описание построения подгруппы $\tildeGamma_{f,0}\subset
\D^0\cap\tildeGamma_f$. } Пусть (после подходящей перенумерации
цилиндров $Z_1,\dots,Z_n$ в~(\ref {eq:Zj}) и соответствующей
перенумерации окружностей $\gamma_1,\dots,\gamma_n$) окружности
$\gamma_\ell\subset M\setminus \N_f \subset M\setminus\N$,
$1\le\ell\le\nu_0$ -- это все такие окружности множества
$\{\gamma_1,\dots,\gamma_n\}$, каждая из которых разбивает
поверхность $M\setminus\N$ на две части (см.\ определение~\ref
{def:F}, обозначение \ref {not:R:G0}, замечание \ref
{rem:numeration}), причем объединение одной из этих двух частей с
окружностью $\gamma_\ell$ гомеоморфно либо кругу, либо проколотому
кругу (т.е.\ кругу без одной внутренней точки), либо цилиндру
$S^1\times[0;1]$. Рассмотрим скручивания Дэна $h_j\in\stab_\D f$
вокруг этих окружностей, $1\le j\le\nu_0$. Каждое такое скручивание
Дэна принадлежит группе $\D^0$, т.е.\ компоненте связности
тождественного диффеоморфизма $\id_M$ в $\Diff(M,\N)$. Значит, все
элементы построенного подмножества
$\{h_1,\dots,h_{\nu_0}\}\subset\{h_1,\dots,h_\ell\}$ принадлежат
группе $\D^0\cap\tildeGamma_f$. Определим подгруппу
$\tildeGamma_{f,0}\subset\tildeGamma_f$ как порожденную
диффеоморфизмами $h_\ell$, $1\le\ell\le\nu_0$.

{\it Описание построения подгрупп
$\tildeGamma_{f,1},\dots,\tildeGamma_{f,e}\subset
\D^0\cap\tildeGamma_f$. } Рассмотрим объединение всех цилиндров в
поверхности $M\setminus\N$, ограниченных парами различных окружностей
из множества $\{\gamma_{\nu_0+1},\dots,\gamma_n\}$ и не содержащих
внутри себя других окружностей этого множества. Это объединение
является либо несвязным объединением $e\ge0$ цилиндров, либо тором (в
этом случае $M=T^2$, $p^*=q^*=r^*=0$ и $\hat p+\hat q+\hat r\ge1$
ввиду (\ref {eq:main}), т.е.\ все критические точки ``плавают'' и по
крайней мере одна из них отмечена, а окружности
$\gamma_{\nu_0+1},\dots,\gamma_n$ попарно изотопны в торе $M$); в
последнем случае положим $e=1$, $\nu_1=n$, и заменим указанное
объединение цилиндров на один цилиндр, содержащий все окружности
$\gamma_{\nu_0+1},\dots,\gamma_n$ и ограниченный двумя из этих
окружностей, причем этот цилиндр не содержит первую отмеченную
критическую точку (эти условия определяют цилиндр однозначно). Пусть
(после подходящей перенумерации цилиндров $Z_{\nu_0+1},\dots,Z_n$
в~(\ref {eq:Zj}) и соответствующей перенумерации окружностей
$\gamma_{\nu_0+1},\dots,\gamma_n$) окружности $\gamma_\ell$,
$\nu_{k-1}<\ell\le\nu_k$ -- это все окружности в $k$-ом из этих $e$
цилиндров, причем можем и будем считать, что нумерация окружностей
идет в порядке следования этих окружностей в $k$-ом цилиндре, $1\le
k\le e$. Композиция $h_\ell^{-1}\circ h_{\ell+1}$ скручиваний Дэна
$h_\ell$ и $h_{\ell+1}$, взятых в противоположных степенях,
принадлежит группе $\D^0$ при $\nu_{k-1}<\ell<\nu_k$, $1\le k\le e$.
При $1\le k\le e$ определим подгруппу
$\tildeGamma_{f,k}\subset\tildeGamma_f$ как порожденную
диффеоморфизмами $h_\ell^{-1}\circ h_{\ell+1}$,
$\nu_{k-1}<\ell<\nu_k$.

Покажем, что группа $(\D^0\cap\tildeGamma_f)^*$ допускает разложение
(\ref {eq:generH}). Осталось показать, что
$\D^0\cap\tildeGamma_f\subseteq\tildeGamma_{f,0}\times\tildeGamma_{f,1}\times\ldots\times\tildeGamma_{f,e}$.
Дополним набор диффеоморфизмов
$h_1,\dots,h_{\nu_0}\in\tildeGamma_{f,0}$, $h_\ell^{-1}\circ
h_{\ell+1}\in\tildeGamma_{f,k}$, $\nu_{k-1}<\ell<\nu_k$, $1\le k\le
e$, до набора образующих группы $\tildeGamma_f\cong\ZZ^n$ набором
скручиваний Дэна $h_i$, $i\in A=A(f)$, где
 $$
 A=A(f):=\{\nu_1,{\nu_2},\dots,{\nu_e},\ {\nu_e+1},{\nu_e+2},\dots,{n}\} \subset
 \{1,\dots,n\}, \quad |A| = n-\nu_e+e.
 $$
Пусть некоторая композиция $h\in\tildeGamma_f$ целых степеней
диффеоморфизмов полученного набора образующих принадлежит группе
$\D^0\cap\tildeGamma_f$. Покажем, что показатель степени каждого из
$|A|=n-\nu_e+e$ диффеоморфизмов $h_i$, $i\in A$, в этой композиции
равен нулю. Произведение $\tilde h$ этих $n-\nu_e+e$ диффеоморфизмов
в тех степенях, в которых они входят в композицию $h$, также является
элементом группы $\D^0\cap\tildeGamma_f$, так как отличается от
исходной композиции $h$ домножением на элемент из подгруппы
$\tildeGamma_{f,0}\times\tildeGamma_{f,1}\times\ldots\times\tildeGamma_{f,e}$,
содержащейся в $\D^0\cap\tildeGamma_f$ по построению. Значит, $\tilde
h\in\D^0\cap\tildeGamma_f\subset\D^0$. Так как окружности
$\gamma_i\subset(\Int M)\setminus\N$, $i\in A$, попарно не
пересекаются, никакая из них не ограничивает цилиндр или (проколотый
или непроколотый) круг в $M\setminus\N$, и никакие две из них не
ограничивают цилиндр в $M\setminus\N$, то скручивания Дэна $h_i$,
$i\in A$, вокруг этих окружностей (рассматриваемые с точностью до
гомотопии в пространстве непрерывных отображений пары $(M,\N)$)
порождают подгруппу группы
$\Homeo(M,\N_0,\N_1,\N_2)/\Homeo^0(M,\N_0,\N_1,\N_2)\cong\D/\D^0$
классов отображений, изоморфную свободной абелевой группе ранга
$|A|=n-\nu_e+e$ (см., например, \cite[лемма 2.1(1)]{BLMC} или
\cite{Pdiplom}). Поэтому показатели степеней всех диффеоморфизмов
$h_i$, $i\in A$, равны 0. Это завершает доказательство разложения
(\ref {eq:generH}).

Построим специальные (криволинейные) координаты в $U_f\subset
U_f^{\infty}$, в которых свободное действие цилиндра
$\RR^n/Z^d\cong\RR^{n-d}\times(S^1)^d$ (см.\ конец шага~6)
``выпрямляется''. Пусть нумерация цилиндров $Z_1,\dots,Z_n$ такая же,
как в~(\ref {eq:generH}). Для любого $u'\in U_f'$ рассмотрим базис
$v_1(u'),\dots,v_n(u')$ в плоскости $\<[\tilde e_1]^*,\dots,[\tilde
e_n]^*\>+u'\subset U_f$ и новый базис
\begin{equation} \label {eq:vu'j}
 \tilde v_i(u'):=v_i(u'), \ i\in A\cup\{1,\dots,\nu_0\}, \quad
 \tilde v_j(u'):=v_j(u')-v_{j+1}(u'),
\ j\in B\setminus\{1,\dots,\nu_0\},
 \end{equation}
а также отвечающее этому базису разложение
$$
 \ker\left[H^1_f\to H^1(G_f,\N_{f,1};\RR)\right]
 = \left\<[\tilde e_\ell]^*\right\>_{\ell=1}^n
 =\left\<\tilde v_{i}(u')\right\>_{i\in A} \oplus \left\<\tilde v_{j}(u')\right\>_{j\in B},
 $$
где $B=B(f):=\{1,\dots,n\}\setminus A$. Тогда для каждого $u'\in
U'_f\subset(U'_f)^\infty$ любой коцикл $\tilde
u=\sum_{\ell=1}^n\tilde u_j[\tilde e_j]^*\in\<[\tilde
e_\ell]^*\>_{\ell=1}^n$ имеет вид
$$
 \tilde u=\sum_{i\in A}x_i \ \tilde v_{i}(u') + \sum_{j\in B}\varphi_j \ \tilde v_{j}(u'),
$$
где координаты $x_i,\varphi_j$ ($i\in A$, $j\in B$) в $n$-мерной
плоскости $\<[\tilde e_1]^*,\dots,[\tilde e_n]^*\>+u'$ выражаются
через координаты $\tilde u_1,\dots,\tilde u_n$,
$u'_{n+1},\dots,u'_{2q}$, по формулам
 \begin{equation}\label{eq:coord}
 x_i=\frac{\tilde u_i}{u'([\gamma_i])}, \quad \nu_e< i\le n, \qquad
 x_{\nu_k}=\frac{\tilde u_{\nu_{k-1}+1}}{u'([\gamma_{\nu_{k-1}+1}])}
 +\ldots+
 \frac{\tilde u_{\nu_k}}{u'([\gamma_{\nu_{k}}])}, \quad 1\le k\le e,
 \end{equation}
$$
 \varphi_j=\frac{\tilde u_j}{u'([\gamma_j])}, \quad 1\le j\le\nu_0, \qquad
 \varphi_j=\frac{\tilde u_{\nu_{k-1}+1}}{u'([\gamma_{\nu_{k-1}+1}])}
 +\ldots+
 \frac{\tilde u_j}{u'([\gamma_{j}])}, \quad
 \nu_{k-1}< j<\nu_k.
$$
Знаменатели в выражениях для $x_i$ и $\varphi_j$ положительны, так
как $u'\in (U'_f)^{\infty}$. Таким образом, на множестве
$U_f^{\infty}$ мы ввели гладкие регулярные координаты $x_i\in\RR$,
$\varphi_j\mod1\in S^1=\RR/\ZZ$ ($i\in A$, $j\in B$),
$(u'_{n+1},\dots,u'_{2q})\in U_f'$. В этих координатах векторные поля
$\tilde v_1,\dots,\tilde v_n$ имеют вид $\tilde v_i=\d/\d x_i$,
$\tilde v_j=\d/\d\varphi_j$ ($i\in A$, $j\in B$), а множества $U_f$ и
$\SS_f$ имеют вид
$$
 U_{[f]_\isot}=U_f \approx (\RR^A \times \RR^B) \times U_f' \approx (\RR^{n-\nu_e+e} \times \RR^{\nu_e-e}) \times U_f' = \RR^n\times U_f', \
$$
$$
 \SS_{[f]_\isot} = \SS_f \approx (\RR^A \times (S^1)^B) \times U_f'
 \approx (\RR^{n-\nu_e+e}\times(S^1)^{\nu_e-e})\times U_f'.
$$
Отсюда действие группы $\tildeGamma_f^*\cong\ZZ^n$ на $U_f$ совпадает
с целочисленными сдвигами вдоль координат $x_i,\varphi_j$ ($i\in A$,
$j\in B$), действие группы
$(\D^0\cap\tildeGamma_f)^*\cong\ZZ^{\nu_e-e}$ на $U_f$ совпадает с
целочисленными сдвигами вдоль координат $\varphi_j$, $j\in B$, а
действие цилиндра
$\RR^A\times(S^1)^B\cong\RR^{n-\nu_e+e}\times(S^1)^{\nu_e-e}$ на
$\SS_f$ совпадает с естественным действием цилиндра сдвигами по себе.
Итак, мы ввели на $\SS_f:=U_f/(\D^0\cap\tildeGamma_f)^*$ структуру
стандартного утолщенного цилиндра (см.\ определение~\ref
{def:thick:cylinder}).

Подмножества $U_f\subset U_f^\infty\subset H^1_f$ инвариантны
относительно правого действия группы
$(\stab_{\D^0}f)^*\subset\Aut(H^1_f)$ на $H^1_f$, а подгруппа
$(\D^0\cap\tildeGamma_f)(\stab_{\D^0}f)^0$ нормальна в
$\stab_{\D^0}f$, где через $(\stab_{\D^0}f)^0$ обозначена подгруппа
группы $\stab_{\D^0}f$, состоящая из всех диффеоморфизмов поверхности
$M$, сохраняющих функцию $f$ и гомотопных $\id_M$ в классе
гомеоморфизмов $M$, сохраняющих функцию $f$. Поэтому имеется
индуцированное правое действие дискретной группы
 \begin{equation} \label{eq:Gamma}
 \Gamma_{[f]_\isot}=\Gamma_f:=(\stab_{\D^0}f)/((\D^0\cap\tildeGamma_f)(\stab_{\D^0}f)^0)
 \end{equation}
на пространствах орбит $\SS_f=U_f/(\D^0\cap\tildeGamma_f)^*$ и
$\SS_f^\infty=U_f^\infty/(\D^0\cap\tildeGamma_f)^*$. Так как действие
группы $\tildeGamma_f$ на $H^0_f$ тривиально (поскольку не
переставляет седловые критические точки), имеем также индуцированное
правое действие группы $\Gamma_f$ на многограннике $D_f$ (см.\ шаг
2).

\begin{Lem} \label{lem:adm}
Если выполнено условие {\rm(\ref {eq:main})}, то индуцированное
покомпонентное правое действие любого диффеоморфизма
$h\in\stab_{\D^0}f$ на прямом произведении
 \begin{equation} \label {eq:stab:adm}
D_f\times\SS_f
 =\tau_{J(\barc(f))}\times\left(U_f/(\D^0\cap\tildeGamma_f)^*\right)
 \approx\tau_{J(\barc(f))}\times\left(\RR^{A(f)}\times(S^1)^{B(f)}\times U_f'\right)
 \end{equation}
является допустимым автоморфизмом стандартной цилиндрической ручки
{\rm(см.\ определение~\ref {def:thick:cylinder}(B))}.
\end{Lem}

\begin{proof}{} Пусть для определенности окружности
$\gamma_\ell\subset Z_\ell$ определены условием
$f(\gamma_\ell)=\frac12(\sup f|_{Z_\ell}+\inf f|_{Z_\ell})$. Тогда
любой диффеоморфизм $h\in\stab_{\D^0}f$ индуцирует перестановку на
множестве окружностей $\gamma_\ell$, $1\le\ell\le n=n(f)$. При этой
перестановке каждая окружность $\gamma_{\nu_0+1},\dots,\gamma_{n}$
переходит в себя (см.\ доказательство леммы~\ref {lem:invar*}, шаг
1). То есть, переставляются только окружности $\gamma_\ell$,
$\ell\in\{1,\dots,\nu_0\}\subset B(f)$ (и отвечающие им векторные
поля $\tilde v_\ell$), а соответствующая перестановка
$\pi\in\Sigma_{|A(f)|}$ тривиальна (см.\ определение~\ref
{def:thick:cylinder}(A)). Если тривиальны также соответствующие
автоморфизмы многогранников $b\:D_f\to D_f$, $a\:U_f'\to U_f'$ и
перестановка $\rho\in\Sigma_{|B(f)|}$, то $h$ переводит в себя каждое
седло, каждое ориентированное ребро графа $G_f$, и каждую окружность
$\gamma_\ell$, а потому принадлежит
$(\D^0\cap\tildeGamma_f)(\stab_{\D^0}f)^0$, откуда его действие на
$D_f\times\SS_f$ совпадает с тождественным отображением. Если
автоморфизм многогранника $b\:D_f\to D_f$ тривиален (что равносильно
тому, что $h$ переводит каждое седло в себя), то перестановка
$\rho\in\Sigma_{|B(f)|}$ окружностей $\gamma_\ell$ тоже тривиальна,
так как в противном случае $h$ нетривиально действует на $H_1(\bar
M)$; тривиальность автоморфизма $a^2\:U'_f\to U'_f$ следует из того,
что $h^2$ переводит каждое ребро графа $G_f$ в себя. Пусть теперь
автоморфизм многогранника $a\:U'_f\to U'_f$ тривиален. Покажем, что
автоморфизм многогранника $b\:D_f\to D_f$ тоже тривиален. Если
количество седловых значений $s(f)>1$, то $U'_f$ является
$(2q-n(f))$-мерным многогранником, поэтому из тривиальности
автоморфизма $a\:U'_f\to U'_f$ следует, что $h$ переводит в себя
каждое ребро графа $G_f$, а потому и каждое седло, откуда $b\:D_f\to
D_f$ тривиален. Если $s(f)=1$, то по лемме \ref {lem:invar*} ниже
отображение $h$ переводит в себя хотя бы одно ребро графа $G_f$ (а
потому и каждое его ребро ввиду связности графа $G_f$, а потому и
каждую седловую точку), откуда автоморфизм $b\:D_f\to D_f$ тривиален.
Так как автоморфизм $b\:D_f\to D_f$ является ограничением
автоморфизма многогранника $\P_f^{q-1}$, индуцированного
перестановкой координатных осей, то по лемме~\ref {lem:P} он
допустим. Лемма~\ref {lem:adm} доказана.
\end{proof}

{\it Шаг 8.} Изучим взаимосвязь утолщенных цилиндров
$\SS_{[f]_\isot},\SS_{[g]_\isot}$ для примыкающих классов изотопности
$[f]_\isot\prec[g]_\isot$. Пусть $f\in F^1$ --- отмеченная функция
Морса класса изотопности $[f]_\isot$. Для любой грани
$\tau'\prec\tau_{J(\barc(f))}=:D_{[f]_\isot}=D_f$ обозначим через
$g\in F^1$ отмеченную функцию класса изотопности
$\delta_{\tau'}[f]_\isot$ (см.~(\ref {eq:incid})) и фиксируем
диффеоморфизм $h_{f,\tau'}\in\D^0$ как в~(\ref {eq:incid}) и~(\ref
{eq:incid'}). Рассмотрим индуцированный изоморфизм
$$
h_{f,\tau'}^*\:H^1_f\to H^1_g
$$
векторных пространств, аналогичный изоморфизму
$h_{f,\tau'}^{*0}\:H^0_f\to H^0_g$ из~(\ref {eq:incid'}). Докажем
включения
\begin{equation}\label{eq:embedU}
 h_{f,\tau'}^*(U_f)\subset U_g, \qquad
 h_{f,\tau'}^*(U_f^{\infty})\subset U_g^{\infty} .
 \end{equation}
Пусть, как и выше, $s=s(f)$ -- количество седловых критических
значений функции $f$, и $k:=q-s$ --- размерность многогранника $D_f$
(см.\ шаги 1, 2). С учетом определения $U_f\subset
U_f^{\infty}\subset H^1_f$ (см.\ (\ref {eq:U}), (\ref {eq:U0})), нам
достаточно показать, что сопряженный к изоморфизму $h_{f,\tau'}^*$
изоморфизм $(h_{f,\tau'})_*\:H_{g,1}=H_1(M\setminus(\N_{g,0}\cup
\N_{g,2}),\N_{g,1};\RR)\to H_1(M\setminus(\N_{f,0}\cup
\N_{f,2}),\N_{f,1};\RR)=H_{f,1}$ переводит гомологический класс
любого ориентированного ребра графа $G_g$ в сумму гомологических
классов некоторых ориентированных ребер графа $G_f$ (см.\ определение
графа $G_f$ в обозначении~\ref{not:Gf:or}), и что количество этих
ребер всегда $\le2k+1$.

Обозначим через $C_g$ компоненту связности графа $G_g$, в которой
лежит рассматриваемое ребро графа $G_g$. Дополнение графа $G_f$ в
поверхности $M$ распадается на ``открытые цилиндры'' $Z_\ell(f)\ccong
S^1\times (0;1)$, $1\le\ell\le n=n([f])$, ``полуоткрытые цилиндры''
$S^1\times [0;1)$ и открытые круги, содержащие ровно одну критическую
точку минимума или максимума функции $f$. Поэтому имеется ретракция
 $$
 \varrho_f\: M_f':=\left(M\setminus(\N_{f,0}\cup \N_{f,2})\right)\setminus\left(\bigcup\limits_{\ell=1}^n\gamma_\ell(f)\right)\to G_f,
 $$
где $\gamma_\ell(f)=S^1\times\{\frac12\}\subset Z_\ell(f)$. Более
точно, определим эту ретракцию так, чтобы она переводила любую точку
поверхности $M_f'$ в точку пересечения проходящей через нее
интегральной траектории векторного поля $\grad f|_{M_f'}$ (в смысле
некоторой фиксированной римановой метрики $ds_0^2$ на $M$) с графом
$G_f$. Без ограничения общности мы также будем считать, что функция
$\tilde f$ и диффеоморфизм $h_{0;f,\tilde f}$ в определении
диффеоморфизма $h_{f,\tau'}$ (см.\ (\ref {eq:incid'})) строились так:
фиксируем попарно непересекающиеся круги вокруг седловых точек
функции $f\in F^1$ и потребуем, чтобы в каждом из них $\tilde
f=f+\const$, и чтобы $h_{0;f,\tilde f}=\id_M$ и функция $\tilde f\in
F^1$ была получена при помощи $C^2$-малого возмущения функции $f$.
Тогда $h_{f,\tau'}(G_g)\subset M_f'$, причем отображение
 \begin{equation} \label{eq:p:fg}
p_{f,\tau'}:=\varrho_f\circ h_{f,\tau'}|_{C_g}\:C_g\to G_f
 \end{equation}
является погружением графов, переводит множество вершин на множество
вершин согласно биекции $h_{f,\tau'}|_{\N_g}\:\N_g\to \N_f$ и
сохраняет ориентацию ребер.
Отсюда получаем, что $p_{f,\tau'}$ переводит любое ориентированное
ребро графа $G_g$ в ориентированный путь на графе $G_f$, ориентация
которого согласована с ориентацией ребер графа $G_f$.

Осталось показать, что длина указанного пути на графе $G_f$ (т.е.\
количество проходимых этим путем ребер графа $G_f$) не превосходит
$2k+1$. Пусть $C_f$ -- компонента связности графа $G_f$, в которой
лежит рассматриваемый путь. Граф $C_f$ имеет не более $k+1$ вершины
(так как число компонент связности графа $G_f$ не меньше чем
$s=q-k$), а потому он имеет не более $2k+2$ ребер. Но наше ребро
является собственным подграфом графа $C_g$, а потому наш путь
является собственным подграфом графа $p_{f,\tau'}(C_g)\subset C_f$.
Так как граф $C_f$ имеет не более $2k+2$ ребер, наш путь имеет не
более $2k+1$ ребер, что и требовалось. Это завершает доказательство
включений (\ref {eq:embedU}).

Изоморфизм $h_{f,\tau'}^*\:
H^1_f\stackrel{\cong}{\longrightarrow}H^1_g$ индуцирует изоморфизм
 $$
\hat
h_{f,\tau'}\:\Aut(H^1_f)\stackrel{\cong}{\longrightarrow}\Aut(H^1_g),
\quad h^*\mapsto h^*_{f,\tau'} h^* (h_{f,\tau'}^*)^{-1}, \qquad
h^*\in\Aut(H^1_f).
 $$
Рассмотрим вложение множеств окружностей
$\{\gamma_\ell(f)\}_{\ell=1}^{n(f)}\hookrightarrow\{\gamma_m(g)\}_{m=1}^{n(g)}$,
сопоставляющее окружности $\gamma_\ell(f)$ окружность
$\gamma_{m(\ell)}(g)$, такую что
$h_{f,\tau'}^{-1}(\gamma_\ell(f))\subset Z_{m(\ell)}(g)$. Тогда для
каждого векторного поля $v_\ell(f)$ на $H^1_f$ (см.~(\ref {eq:v:ell})
и~(\ref {eq:v})) выполнено
$(h_{f,\tau'}^*)_*(v_\ell(f))=v_{m(\ell)}(g)$. Отсюда $\hat
h_{f,\tau'}(\tildeGamma_f^*)\subset\tildeGamma_g^*$, а потому
 $$
 \hat h_{f,\tau'}((\D^0\cap\tildeGamma_f)^*)\subset(\D^0\cap\tildeGamma_g)^*.
 $$
Поэтому вложение $h_{f,\tau'}^*|_{U_f}\: U_f\hookrightarrow U_g$
(см.~(\ref {eq:embedU})) индуцирует корректно определенное
отображение пространств орбит
\begin{equation}\label{eq:chi}
 [h_{f,\tau'}^*|_{U_f}]\:U_f/(\D^0\cap\tildeGamma_f)^*\looparrowright U_g/(\D^0\cap\tildeGamma_g)^*,
\end{equation}
являющееся погружением утолщенных цилиндров (так как группа
$\tildeGamma_f^*$ действует свободно и дискретно на $U_f$, см.\ шаги
6 и 7).

Докажем, что погружение (\ref {eq:chi}) утолщенных цилиндров является
допустимым (см.\ определение~\ref {def:thick:cylinder},(D)). Из
$(h_{f,\tau'}^*)_*(v_\ell(f))=v_{m(\ell)}(g)$, $1\le\ell\le n$, (\ref
{eq:vu'j}) и описания подгруппы
$\tildeGamma_{f,k}\subset\D^0\cap\tildeGamma_f$ (см.\ шаг 7) следует,
что при $\nu_{k-1}(f)<j<\nu_k(f)$, $1\le k\le e(f)$, выполнено
$$
 (h_{f,\tau'}^*)_*(\tilde v_j(f))
 =(h_{f,\tau'}^*)_*(v_j(f)-v_{j+1}(f))
 = v_{m(j)}(g)-v_{m(j+1)}(g)
$$
$$
=(v_{m(j)}(g)-v_{m(j)+\eta_k}(g))+(v_{m(j)+\eta_k}(g)-v_{m(j)+2\eta_k}(g))+\ldots
$$
$$
 =\eta_k(\tilde v_{m(j)+(\eta_k-1)/2}(g)+\tilde v_{m(j)+(3\eta_k-1)/2}(g)+\ldots+\tilde v_{m(j+1)-(1+\eta_k)/2}(g)),
$$
где $\eta_k=\eta_k(f,g):=\sgn(m(\nu_k-1)-m(\nu_k))$. Отсюда следует,
что при вложении $h_{f,\tau'}^*|_{U_f}\: U_f\hookrightarrow U_g$
коммутирующие векторные поля (a) $\tilde v_i(f)$, $i\in A(f)$, (b)
$\tilde v_j(f)$, $j\in B(f)$, на $U_f$ (потоки которых задают
свободное действие группы $\RR^{A(f)}\times\RR^{B(f)}$ на $U_f$ и
свободное действие цилиндра $\RR^{A(f)}\times(S^1)^{B(f)}$ на
$U_f/(\D^0\cap\tildeGamma_f)^*$) переходят в следующие векторные поля
на $U_g$:

(a) $\tilde v_{m(i)}(g)$ (при $\nu_e(g)<m(i)\le n(g)$) или $\tilde
v_{m(i)}(g)+\tilde v_{m(i)+1}(g)+\ldots+\tilde v_{\nu_t(g)}$ (при
$\nu_{t-1}(g)<m(i)\le\nu_t(g)$, $1\le t\le e(g)$),

(b) $\tilde v_{m(j)}(g)$ (при $1\le j\le\nu_0(f)$) или $\eta_k(\tilde
v_{m(j)+(\eta_k-1)/2}(g)+\tilde
v_{m(j)+(\eta_k-1)/2+\eta_k}(g)+\ldots+\tilde
v_{m(j+1)-(1+\eta_k)/2}(g))$ (при $\nu_{k-1}(f)<j<\nu_k(f)$, $1\le
k\le e(f)$),
\\
причем каждое поле $\tilde v_m(g)$, $1\le m\le n(g)$, входит в
качестве слагаемого (с коэффициентом $\pm$) не более чем в одно из
полей $(h_{f,\tau'}^*)_*(\tilde v_\ell(f))$, $1\le\ell\le n(f)$.

Из описанного поведения векторных полей $\tilde v_\ell$, $1\le\ell\le
n$, при погружении~(\ref {eq:chi}) следует, что это погружение
является допустимым погружением утолщенных цилиндров (см.\
определение~\ref {def:thick:cylinder}(D)).

\subsection{Построение косой цилиндрической ручки $\DD_{[f]_\isot}^\st$ для класса изотопности $[f]_\isot$}\label{subsec:handle}

{\it Шаг 9.} Рассмотрим стандартную цилиндрическую ручку
$D_{[f]_\isot}\times\SS_{[f]_\isot}=D_f\times\SS_f=\tau_{J(\barc(f))}\times(U_f/(\D^0\cap\tildeGamma_f)^*)$
и покомпонентное правое действие на ней дискретной группы
 $$
 \Gamma_{[f]_\isot}=\Gamma_f\cong\widetilde\Gamma_f /\left(((\D^0\cap\tildeGamma_f)(\stab_{\D^0}f)^0)/(\stab_{\D^0}f)^0\right),
 $$
где
\begin{equation}\label{eq:Gamma:f}
 \widetilde\Gamma_{[f]_\isot}=\widetilde\Gamma_f:=(\stab_{\D^0}f)/(\stab_{\D^0}f)^0,
 \end{equation}
допустимыми автоморфизмами (см.\ (\ref {eq:Gamma}) и лемму~\ref
{lem:adm}).
Покажем, что это действие (а также действие группы
$\widetilde\Gamma_f$ на $U_f^\infty$) свободно. Докажем две леммы.

\begin{Lem} \label{lem:invar*}
Если выполнено условие~{\rm(\ref {eq:main})}, то для любого
диффеоморфизма $h\in\stab_{\T}f$ найдется ребро графа $G_f$ {\rm(см.\
обозначения~\ref {not:R:G0}(B) и~\ref {not:Gf:or})}, переходящее в
себя при отображении $h$.
\end{Lem}

\begin{proof}{}
{\it Шаг 1.} Пусть $W_f$ -- граф Кронрода-Риба функции $f$
(см.~\cite{Kronrod} или~\cite {Kmsb}), т.е.\ граф $W_f$ получен из
поверхности $M$ стягиванием в точку каждой компоненты связности линии
уровня функции $f$. Обозначим через $p_f\:M\to W_f$ естественную
проекцию. Вершину графа $W_f$ назовем {\it сферической}, если прообраз
достаточно малой ее окрестности при отображении $p_f$ гомеоморфен
сфере с проколами. Вершину графа $W_f$ назовем {\it граничной}
(соответственно {\it отмеченной}), если ее прообраз при отображении
$p_f$ является компонентой края $M$ (соответственно содержит
отмеченную критическую точку функции $f$). Подграф графа $W_f$
назовем {\it $\stab_{\T}f$-неподвижным}, если при автоморфизме
$p_f\circ h\circ p_f^{-1}$ графа $W_f$, индуцированном любым
диффеоморфизмом $h\in\stab_{\T}f$, любая вершина и любое ребро этого
подграфа переходят в себя.
Обозначим через $W_f'$ минимальный связный подграф графа $W_f$,
содержащий каждый простой цикл графа $W_f$, каждую граничную вершину,
каждую отмеченную вершину и каждую несферическую вершину. Он непуст и
содержит неграничную вершину в силу~(\ref {eq:main}). Покажем, что
подграф $W_f'$ является $\stab_{\T}f$-неподвижным.
Пусть $h\in\stab_{\T}f$. Так как $h$ сохраняет неподвижными все
отмеченные критические точки функции $f$ и переводит в себя каждую
компоненту края, то все отмеченные вершины и все граничные вершины
$\stab_{\T}f$-неподвижны. Ввиду $h\in\T$ индуцированный автоморфизм
гомологий $\bar h^*\in H_1(\bar M)$ (см.\ обозначение~\ref
{not:R:G0}(B)) совпадает с тождественным, поэтому каждая
несферическая вершина $v\in W_f$ является $\stab_{\T}f$-неподвижной,
а каждый простой цикл на графе $W_f$ переходит в себя с сохранением
ориентации при отображении $p_f\circ h\circ p_f^{-1}$. Если
пересечение двух простых циклов непусто и связно, то оно
$\stab_{\T}f$-неподвижно, поэтому такие циклы
$\stab_{\T}f$-неподвижны. Поэтому каждая компонента связности
объединения простых циклов, не являющаяся простым циклом (или
содержащая несферическую или отмеченную вершину),
$\stab_{\T}f$-неподвижна. Пусть $W_f''\subset W_f$
--- объединение всех простых циклов, множества несферических вершин и
множества отмеченных вершин графа $W_f$, и пусть простой путь в графе
$W_f$ соединяет две компоненты связности графа $W_f''$ и пересекается
с $W_f''$ только в концах. Такой путь единствен (для фиксированной
пары компонент), а потому $\stab_{\T}f$-неподвижен (в силу $p_f\circ
h\circ p_f^{-1}$-инвариантности подграфа $W_f''$). Отсюда следует $\stab_{\T}f$-неподвижность
подграфа $W_f'$.

{\it Шаг 2.} Пусть $M'\subset M$ --- прообраз малой связной
окрестности вершины $v\in W_f$ при отображении $p_f\:M\to W_f$,
обозначим $c:=f(p_f^{-1}(v))$. Пусть $h\in\stab_\T f$, $h(M')=M'$, и
существует седловая критическая точка в $M'$, не являющаяся
неподвижной при $h$.
Пусть замкнутая поверхность $\overline{M'}$ получена из $M'$
стягиванием в точку каждой компоненты края $M'$, и $\bar
h\:\overline{M'}\to\overline{M'}$ --- индуцированный гомеоморфизм.
Без ограничения общности будем считать, что ограничение $h$ на каждую
$h$-инвариантную компоненту связности $\overline{M'}\setminus
f^{-1}(c)$, гомеоморфную
$(S^1\times(c;c+\eps])/(S^1\times\{c+\eps\})$, является поворотом
вокруг точки $S^1\times\{c+\eps\}$ (по отношению к естественным
``полярным'' координатам в
$(S^1\times(c;c+\eps])/(S^1\times\{c+\eps\})$), тогда эта точка
является единственной неподвижной точкой в данной компоненте. Так как
$\bar h$ индуцирует тождественный автоморфизм гомологий
$H_1(\overline{M'})$ (ввиду $h\in\T$), и все его неподвижные точки
имеют индекс $+1$, то по формуле Лефшеца количество неподвижных точек
равно $\chi(\overline{M'})$.

{\it Шаг 3.} Вершину графа $W_f$ назовем {\it сильно
$\stab_{\T}f$-инвариантной}, если ее степень в графе $W_f$ больше 1 и
при любом диффеоморфизме $h\in\stab_{\T}f$ каждая вершина и каждое
ребро графа $p_f^{-1}(v)\subset G_f$ переходят в себя. Докажем, что в
$W_f'$ существует сильно $\stab_{\T}f$-инвариантная вершина. Если
вершина $v\in W_f'$ не является сильно $\stab_{\T}f$-инвариантной, то
либо ее степень в $W_f$ равна 1, либо найдется такой диффеоморфизм
$h\in\stab_\T f$, что количество неподвижных точек соответствующего
индуцированного гомеоморфизма $\bar h\:\overline{M'}\to\overline{M'}$
(см.\ шаг 2 выше) согласно формуле Лефшеца равно
$\chi(\overline{M'})$ и не меньше суммы $k_v+\deg_{W_f'}v$ числа
$k_v$ отмеченных критических точек в $p_f^{-1}(v)$ и степени
$\deg_{W_f'}v$ вершины $v$ в графе $W_f'$ (т.е.\
$\chi(\overline{M'})\ge k_v+\deg_{W_f'}v\ge0$), откуда вершина $v$
является сферической и $k_v+\deg_{W_f'}v\le2$ (так как в противном
случае $\chi(\overline{M'})=k_v+\deg_{W_f'}v=0$, откуда $W_f'=\{v\}$,
$M$ --- тор и все критические точки неотмечены, что противоречит
(\ref {eq:main})). Если все вершины графа $W_f'$ не являются сильно
$\stab_{\T}f$-инвариантными, то по доказанному выше каждая вершина
$v\in W_f'$ сферическая и либо имеет степень 1 в $W_f$, либо имеет
степень $\deg_{W_f'}v\le2-k_v\le2$ в $W_f'$, откуда граф $W_f'$
является простой ломаной, все его внутренние вершины неотмечены (так
как $k_v=0$ в случае $\deg_{W_f'}v=2$), а сумма значений $k_v$ для
концевых вершин $v\in\d W_f'$, не являющихся граничными, не
превосходит $2-d^+-d^-=\chi(M)$ (так как $k_v\le1$ при
$\deg_{W_f'}v=1$, и $k_v\le2$ при $\deg_{W_f'}v=0$), откуда общее
количество отмеченных критических точек $\le\chi(M)$, что
противоречит (\ref {eq:main}). Лемма~\ref {lem:invar*} доказана.
\end{proof}

\begin{Lem} \label{lem:free}
Пусть выполнено условие~{\rm(\ref {eq:main})}. Пусть функция Морса
$f\in F^1$, класс относительных 1-когомологий $u\in U_f^{\infty}$,
диффеоморфизм $h\in\stab_{\T}f$ {\rm(см.\ обозначение~\ref
{not:R:G0}(B))} и мультискручивание Дэна $h_1\in\tildeGamma_f$
удовлетворяют условию $h^*(u)=h_1^*(u)$. Тогда
$hh_1^{-1}\in(\stab_{\D^0}f)^0$, см.~{\rm(\ref {eq:Gamma})}.
\end{Lem}

\begin{proof}{}
По лемме~\ref {lem:invar*} найдется ребро графа $G_f$, переходящее в
себя при отображении $h$. Пусть $C_f$ --- компонента связности графа
$G_f$, содержащая это ребро, и пусть $Z_\ell$ --- открытый цилиндр,
одна из компонент границы которого (скажем, нижнее основание
$\d^-Z_\ell$) имеет общее ребро с графом $C_f$ (см.~(\ref {eq:Zj})).
Тогда любое ребро графа $C_f\cup\d^-Z_\ell$ и цилиндр $Z_\ell$ тоже
переходят в себя при отображении $h$. Так как пути $\tilde
e_\ell,h_1(\tilde e_\ell)\subset\overline{Z_\ell}$ выходят из одной и
той же точки (принадлежащей $\d^-Z_\ell$), то 1-цепь $h(\tilde
e_\ell)-h_1(\tilde e_\ell)$ гомологична некоторой линейной комбинации
$\sum_{i=1}^{2q}\lam_ie_i$ ориентированных ребер основания
$\d^+Z_\ell$ с целыми коэффициентами, причем все коэффициенты
$\lam_i$ либо неотрицательны, либо неположительны одновременно. Но
 $$
 \sum_{i=1}^{2q}\lam_iu([e_i])= u([h(\tilde e_\ell)]-[h_1(\tilde e_\ell)])
 =(h^*(u)-h_1^*(u))([\tilde e_\ell])=0
 $$
по предположению. Так как значение 1-коцикла $u\in U_f^{\infty}$ на
каждом ориентированном ребре \\ $e_1,\dots,e_{2q}\subset G_f$
положительно (см.~(\ref {eq:U0})), то линейная комбинация тривиальна.
Значит, $[h(\tilde e_\ell)]=[h_1(\tilde e_\ell)]$, откуда
$hh_1^{-1}|_{\overline{Z_\ell}}$ гомотопно $\id_{\overline{Z_\ell}}$
в классе гомеоморфизмов, сохраняющих функцию $f|_{\overline{Z_\ell}}$
и переводящих вершины графа $\d Z_\ell$ в себя. Эти рассуждения
показывают (с использованием связности $M$), что
$hh_1^{-1}\in(\stab_{\D^0}f)^0$. Лемма доказана.
\end{proof}

В силу (\ref {eq:Gamma}) и леммы \ref {lem:free} группа $\Gamma_f$
действует свободно на утолщенном цилиндре
$\SS_f=U_f/(\D^0\cap\tildeGamma_f)^*$, поэтому она действует свободно
на стандартной цилиндрической ручке $D_f\times\SS_f$ допустимыми
автоморфизмами (см.\ лемму~\ref {lem:adm}), а потому конечна (см.\
определение \ref {def:thick:cylinder}(B)). Значит, пространство орбит
 \begin{equation} \label{eq:DD}
 \DD_{[f]_\isot}^\st=\DD_f^\st:=
 (D_f\times\SS_f)/\Gamma_f
 \approx (\tau_{J(\barc(f))}\times U_f)/\widetilde\Gamma_f
 \end{equation}
является стандартной косой цилиндрической ручкой (см.~(\ref
{eq:Gamma:f}) и определение~\ref {def:thick:cylinder}(C)).

{\it Шаг 10.} Изучим взаимосвязь стандартных косых цилиндрических
ручек $\DD^\st_{[f]_\isot},\DD^\st_{[g]_\isot}$ для примыкающих
классов изотопности $[f]_\isot\prec[g]_\isot$. Пусть $f\in F^1$ ---
отмеченная функция Морса класса изотопности $[f]_\isot$. Для любой
грани $\tau'\prec\tau_{J(\barc(f))}=D_{[f]_\isot}=D_f$ обозначим
через $g\in F^1$ отмеченную функцию класса изотопности
$[g]_\isot=\delta_{\tau'}[f]_\isot$ (см.~(\ref {eq:incid})).
Рассмотрим погружение
 \begin{equation} \label {eq:face}
 h_{f,\tau'}^{*0}|_{\tau'}\times[h_{f,\tau'}^*|_{U_f}]\: \tau'\times\SS_f \looparrowright D_g\times\SS_g,
 \end{equation}
 $$
 \left(\barc,(\D^0\cap\tildeGamma_f)^*(u)\right) \mapsto \left(h_{f,\tau'}^{*0}(\barc),(\D^0\cap\tildeGamma_g)^*h_{f,\tau'}^*(u)\right),
 $$
являющееся прямым произведением изометрии~(\ref {eq:glue:D})
евклидовых многогранников и допустимого погружения~(\ref {eq:chi})
утолщенных цилиндров, т.е.\ допустимым погружением стандартных
цилиндрических ручек (см.\ определение~\ref {def:thick:cylinder}(D)).
Рассмотрим орбиту грани $\tau'\subset\d D_f $ при действии группы
$\Gamma_f$ (см.~(\ref {eq:Gamma:f})) изометриями многогранника $D_f$,
и следующие два объединения его граней:
 \begin{equation} \label {eq:delta}
\Gamma_f(\tau'):=\bigcup_{h\in\stab_{\D^0}f}h^{*0}(\tau')
 \ \subseteq \ \bigcup_{\delta_{\tau'_1}[f]_\isot=[g]_\isot
 } \tau'_1 \ =: \ \d_{[g]_\isot}D_{[f]_\isot}=\d_gD_f,
 \end{equation}
см.~(\ref {eq:incid}). Включение в (\ref {eq:delta}) следует из того,
что $\delta_{h^{*0}(\tau')}[f]_\isot=[gh]_\isot=[g]_\isot$ ввиду
включения $h\in\D^0$. Итак, (допустимые) погружения, отвечающие этим
граням, имеют одну и ту же область значений -- стандартную
цилиндрическую ручку $D_g\times\SS_g$. Любые две такие грани либо
совпадают, либо не пересекаются в силу леммы~\ref {lem:P}. Рассмотрим
погружение, составленное из (допустимых) погружений этих граней:
 $$
\left(\d_{g}D_f\right)\times\SS_f
 \looparrowright D_g\times\SS_g, \quad
 \left(\barc,(\D^0\cap\tildeGamma_f)^*(u)\right)
 \mapsto
  \left(h_{f,\tau'_1}^{*0}(\barc),(\D^0\cap\tildeGamma_g)^*h_{f,\tau'_1}^*(u)\right),
 $$
$(\barc,u)\in\tau'_1 \times U_f$. Это отображение корректно
определено (и является погружением), так как грани
$\tau'_1\in(\delta[f]_\isot)^{-1}([g]_\isot)$ попарно не пересекаются
(см.\ выше). Оно переводит любую $\Gamma_f$-орбиту в некоторую
$\Gamma_g$-орбиту, так как ввиду~(\ref {eq:stab}) для любого
$h\in\stab_{\D^0}f$ точки $(\barc,u)\in\tau'_1\times U_f$ и
$(h^{*0}(\barc),h^*(u))\in(h^{*0}(\tau'_1))\times U_f$ переходят в
элементы
         $(h_{f,\tau'_1}^{*0}(\barc),(\D^0\cap\tildeGamma_g)^*h_{f,\tau'_1}^*(u))$
и \\
 $(h_{f,h^{*0}(\tau'_1)}^{*0}h^{*0}(\barc),(\D^0\cap\tildeGamma_g)^*h_{f,h^{*0}(\tau'_1)}^*h^*(u))
 =(h_1^{*0}h_{f,\tau'_1}^{*0}(\barc),(\D^0\cap\tildeGamma_g)^*h_1^*h_{f,\tau'_1}^*(u))$
одной и той же $\Gamma_g$-ор\-би\-ты, где $h_1\in\stab_{\D^0}g$. Поэтому
это погружение индуцирует корректно определенное погружение
пространств орбит:
 $$
  \chi_{[f]_\isot,[g]_\isot}= \chi_{f,g}
 \:
 \left(\left(\d_gD_f\right)\times\SS_f\right)/\Gamma_f
 \looparrowright
(D_g\times\SS_g)/\Gamma_g,
 $$
 $$
 \Gamma_f(\barc,(\D^0\cap\tildeGamma_f)^*(u)) \mapsto \Gamma_g(h_{f,\tau'_1}^{*0}(\barc),(\D^0\cap\tildeGamma_g)^*h_{f,\tau'_1}^*(u)),
 \qquad (\barc,u)\in\tau'_1\times U_f,
 $$
где $[g]_\isot=\delta_{\tau'_1}[f]_\isot$ (см.~(\ref {eq:incid})).
(Оно является погружением, так как группы $\Gamma_f,\Gamma_g$ конечны
и действуют свободно.) Рассмотрим его ограничение:
 $$
 \chi_{[f]_\isot,\tau'}= \chi_{f,\tau'} := \chi_{f,g}|_{\left(\left(\Gamma_f(\tau')\right)\times\SS_f\right)/\Gamma_f}
 \:
 \left(\left(\Gamma_f(\tau')\right)\times\SS_f\right)/\Gamma_f
 \stackrel{(*)}{\approx} \left( \tau' \times\SS_f \right) / \Gamma_{f,\tau'}
 \approx
 $$
 $$
 \approx \left( \tau' \times U_f \right) / \widetilde\Gamma_{f,\tau'}
 \looparrowright
 \DD^\st_{[g]_\isot}=\DD^\st_g=(D_g\times\SS_g)/\Gamma_g
 \approx (D_g\times U_g)/\widetilde\Gamma_g,
 $$
 $$
\widetilde\Gamma_{f,\tau'}(\barc,u) \mapsto
\widetilde\Gamma_g(h_{f,\tau'}^{*0}(\barc),h_{f,\tau'}^*(u)), \quad
(\barc,u)\in\tau'\times U_f,
$$
$$
\Gamma_{f,\tau'}:=
\widetilde\Gamma_{f,\tau'}/\left(((\D^0\cap\tildeGamma_f)(\stab_{\D^0}f)^0)/(\stab_{\D^0}f)^0\right),
\quad \widetilde\Gamma_{f,\tau'}:=\stab_{\widetilde\Gamma_f}\tau',
$$
где гомеоморфизм $(*)$ следует из того, что при $h\in\stab_{\D^0}f$
грани $\tau',h^{*0}(\tau')$ либо совпадают, либо не пересекаются
(см.\ выше). Итак, областью определения погружения $\chi_{f,\tau'}$
является косая грань
\begin{equation} \label {eq:delta:DD}
 \d_{\tau'}\DD^\st_{[f]_\isot}=\d_{\tau'}\DD^\st_f := \left(\left(\Gamma_f(\tau')\right)\times\SS_f\right)/\Gamma_f
\end{equation}
стандартной косой цилиндрической ручки
$\DD^\st_f=(D_f\times\SS_f)/\Gamma_f$, т.е.\ образ грани
$\tau'\times\SS_f$ стандартной цилиндрической ручки $D_f\times\SS_f$
при проекции $D_f\times\SS_f\to\DD^\st_f$. А областью определения
погружения $\chi_{f,g}$ является объединение попарно непересекающихся
косых граней ручки $\DD^\st_f$:
\begin{equation} \label {eq:delta:g:DD}
 \d_{[g]_\isot}\DD^\st_{[f]_\isot}=\d_{g}\DD^\st_f := \left(\left(\d_gD_f\right)\times\SS_f\right)/\Gamma_f.
\end{equation}
Подчеркнем, что
$\d_{\tau'}\DD^\st_{[f]_\isot}=\d_{\tau'_1}\DD^\st_{[f]_\isot}$,
$\chi_{f,\tau'}=\chi_{f,\tau'_1}$ для любой грани
$\tau'_1\in\Gamma_f(\tau')$.

Пусть $[f]_\isot\prec[g]_\isot\prec[g_1]_\isot$, причем
$[g]_\isot=\delta_{\tau'}[f]_\isot$,
$[g_1]_\isot=\delta_{\tau''}[f]_\isot$ для некоторых граней
$\tau''\prec\tau'\prec\tau_{J(\barc)}$, и пусть $f,g,g_1$ ---
отмеченные функции своих классов изотопности. Из~(\ref
{eq:obstruction}) и того, что группы $\Diff^0(M,\N_{g_1})$ и
$\stab_{\D^0}g_1$ действуют тривиально на косой ручке
$\DD^\st_{g_1}$, получаем
 \begin{equation} \label {eq:obstruction:no}
  \chi_{g,g_1} \circ \chi_{f,g} |_{\d_{\tau''}\DD^\st_f} =
 \chi_{g,h_{f,\tau'}^{*0}(\tau'')} \circ \chi_{f,\tau'} |_{\d_{\tau''}\DD^\st_f} = \chi_{f,\tau''}
 = \chi_{f,g_1} |_{\d_{\tau''}\DD^\st_f} .
 \end{equation}

Покажем, что погружение $\chi_{f,g}$ является вложением (а потому
$\chi_{f,\tau'}$ является мономорфизмом стандартных косых
цилиндрических ручек, см.\ определение~\ref {def:thick:cylinder}(D),
ввиду допустимости погружения (\ref {eq:face})). Предположим, что
$u_1,u_2\in U_f^\infty$ и
$h_{f,\tau'}^*(u_1)=h^*h_{f,\tau'_1}^*(u_2)$ для некоторых
$h\in\stab_{\D^0}g$ и $\tau'_1\prec D_f$, таких что
$\delta_{\tau'_1}[f]_\isot=\delta_{\tau'}[f]_\isot=[g]_\isot$.
Покажем, что $u_1\in\widetilde\Gamma_f^*(u_2)$, где
$\widetilde\Gamma_f^*\subset\Aut(H^1_f)$ --- группа автоморфизмов
относительных когомологий, индуцированная группой классов отображений
$\widetilde\Gamma_f=(\stab_{\D^0}f)/(\stab_{\D^0}f)^0$. Имеем
$$
u_1=(h_{f,\tau'}^*)^{-1}h^*h_{f,\tau'_1}^*(u_2)=(h_{f,\tau'_1}hh_{f,\tau'}^{-1})^*(u_2)
$$
\begin{equation}\label{eq:u}
 =(h_{0;f,\tilde f_1}h_{1;\tilde f_1,g}hh_{1;\tilde f,g}^{-1}h_{0;f,\tilde f}^{-1})^*(u_2)
 =(h_{0;f,\tilde f_1}h_{1;\tilde f_1,\tilde f}h_{0;f,\tilde f}^{-1})^*(u_2)
 =h_1^*(u_2),
\end{equation}
где ``возмущенным'' функциям $\tilde f,\tilde f_1$ отвечают грани
$\tau',\tau'_1$ по правилу~(\ref {eq:incid}), $h_{1;\tilde f_1,\tilde
f}:=h_{1;\tilde f_1,g}hh_{1;\tilde f,g}^{-1}$, $h_1:=h_{0;f,\tilde
f_1}h_{1;\tilde f_1,\tilde f}h_{0;f,\tilde f}^{-1}$. Так как
диффеоморфизм $h_{1;\tilde f_1,\tilde f}\in\D^0$ переводит линии
уровня функции $\tilde f$ в линии уровня функции $\tilde f_1$ с
сохранением направления роста (в силу $h\in\stab_{\D^0}g$), то
$\tilde f=h_2\tilde f_1h_{1;\tilde f_1,\tilde f}$ для некоторого
$h_2\in\Diff^+[-1;1]$, откуда
 $\tilde f
 =h_2\tilde f_1(h_{0;f,\tilde f_1}^{-1}h_1h_{0;f,\tilde f})$,
т.е.\ диффеоморфизм $h_1\in\D^0$ переводит линии уровня функции
$\tilde f^*:=\tilde fh_{0;f,\tilde f}^{-1}$ в линии уровня функции
$\tilde f_1^*:=\tilde f_1h_{0;f,\tilde f_1}^{-1}$ с сохранением
направления роста. Отсюда и из того, что обе ``возмущенные'' функции
$\tilde f^*,\tilde f_1^*$ близки к $f$ и имеют те же критические
точки, что и ``невозмущенная'' функция $f$, следует, что
$h_1\in\Diff(M,\N_{f,0},\N_{f,1},\N_{f,2})$. Осталось доказать, что
$h_1\in(\stab_{\D^0}f)(\Diff^0(M,\N_f))$. Так как $\tilde
f^*=h_2\tilde f_1^*h_1$, то
\begin{equation}\label{eq:tilde}
 h_1(G_{\tilde f^*})=G_{\tilde f_1^*}, \ J(\barc(\tilde f_1^*h_1))=J(\barc(\tilde f^*))=:\hat J,
 \quad \Rightarrow
 \
J(\barc(fh_1))=J(\barc(f))=:J
\end{equation}
(так как из того, что перестановка $h_1|_{\N_{f,1}}\:\N_{f,1}\to
\N_{f,1}$ седловых точек переводит отношение частичного порядка
$J(\barc(\tilde f^*_1))$ на множестве $\N_{f,1}$ значениями одной
``возмущенной'' 0-коцепи $\barc(\tilde f^*_1)=\tilde
f^*_1|_{\N_{f,1}}\in C^0(\N_{f,1};\RR)$ в отношение частичного
порядка $J(\barc(\tilde f^*))$ на множестве $\N_{f,1}$ значениями
другой ``возмущенной'' 0-коцепи $\barc(\tilde f^*)=\tilde
f^*|_{\N_{f,1}}\in C^0(\N_{f,1};\RR)$, следует сохранение этой
перестановкой ``более слабого'' отношения частичного порядка на
$\N_{f,1}$ значениями ``невозмущенной'' 0-коцепи
$\barc(f)=f|_{\N_{f,1}}\in C^0(\N_{f,1};\RR)$). Так как функции
$f,fh_1$ имеют одни и те же множества $\N_{f,0}$, $\N_{f,1}$,
$\N_{f,2}$ критических точек минимумов, седловых точек и точек
максимумов, то, согласно (\ref {eq:tilde}) и достаточному условию
изотопности функций Морса (см.~\cite[лемма 1]{KP2}), достаточно
показать совпадение графов $G_{fh_1}=h_0(G_f)$ для некоторого
диффеоморфизма $h_0\in\Diff^0(M,\N_f)$, см.\ обозначение~\ref
{not:Gf:or} и~(\ref {eq:obstruction}) (т.е.\ что графы $G_f$ и
$G_{fh_1}$ изотопны в поверхности $M':=M\setminus(\N_{f,0}\cup
\N_{f,2})$ относительно множества вершин $\N_{f,1}$). Обозначим
$\N_{f,0,2}:=\N_{f,0}\cup \N_{f,2}$.

\begin{Lem} \label{lem:tilde}
Пусть $f\in F^1$, $u\in U_f^\infty\subset H^1_f$, и пусть в
обозначениях~{\rm(\ref {eq:incid})} возмущенной функции $\tilde
f^*:=\tilde fh_{0;f,\tilde f}^{-1}$ отвечает грань
$\tau':=\tau_{J(\barc(\tilde f^*))}\prec\tau_{J(\barc(f))}$. Тогда
граф $G_f\subset M$ {\rm(см.\ обозначение~\ref {not:Gf:or})}
совпадает с точностью до диффеоморфизмов из $\Diff^0(M,\N_f)$ с
некоторым графом $G:=G_{M,\N_{f,0,2},G_{\tilde f^*},u,\hat
J,J}\subset M$, рассматриваемым с точностью до диффеоморфизмов из
$\Diff^0(M,\N_{f,0,2}\cup V(G))$ и определяемым следующим набором
данных: (i) поверхность $M$; (ii) подмножество
$\N_{f,0,2}:=\N_{f,0}\cup \N_{f,2}\subset M$; (iii) граф $G_{\tilde
f^*}\subset M':=M\setminus \N_{f,0,2}$, рассматриваемый с точностью
до диффеоморфизмов из $\Diff^0(M,\N_{f,0,2}\cup V(G_{\tilde f^*}))$;
(iv) класс относительных 1-когомологий $u\in U_f^\infty\subset
H^1_f=H^1(M',V(G_{\tilde f^*});\RR)$; (v) два отношения частичного
порядка $\hat J:=J(\barc(\tilde f^*))\prec J:=J(\barc(f))$ на
множестве $V(G_{\tilde f^*})$ вершин графа $G_{\tilde f^*}$
значениями функций $f,\tilde f^*$ на вершинах (соответственно), где
через $V(G)\subset G$ обозначено множество вершин графа $G$. То есть,
если набор данных (i)--(v) построен указанным способом по паре
``невозмущенной'' и ``возмущенной'' функций $f,\tilde f^*$, имеющих
одно и то же множество критических точек, то
$G_f\in(\Diff^0(M,\N_f))(G)$ для $G:=G_{M,\N_{f,0,2},G_{\tilde
f^*},u,\hat J,J}$ (т.е.\ графы $G_f$ и $G$ изотопны в поверхности
$M'$ относительно множества вершин $V(G_{\tilde f^*})=V(G)$).
\end{Lem}

\begin{proof}{} Сначала проведем доказательство в случае, когда
$\tau':=\tau_{\hat J}$ является гипергранью грани $\tau:=\tau_J$.
Обозначим ``возмущенную'' функцию через $\tilde f_1:=\tilde f^*$.
Граф $G=G_{M,\N_{f,0,2},G_{\tilde f_1},u,\hat J,J}\subset M$ строится
так. Так как $\tau':=\tau_{\hat J}\prec\tau:=\tau_J$ -- гипергрань,
то ровно одно из седловых критических значений $c\in f(\N_{f,1})$
``невозмущенной'' функции $f$ распадается на два седловых критических
значения $c^-<c^+$ ``возмущенной'' функции $\tilde f_1$. Пусть
$Z_\ell\subset M\setminus G_{\tilde f_1}$ -- такая компонента
связности множества $M\setminus G_{\tilde f_1}$, что $\inf\tilde
f_1|_{Z_\ell}=c^-$ и $\sup\tilde f_1|_{Z_\ell}=c^+$ (т.е.\ $Z_\ell$
является открытым цилиндром для функции $\tilde f_1$, см.~(\ref
{eq:Zj})). Ориентированный граф $G\subset\overline{Z_\ell}$ назовем
{\it $Z_\ell$-допустимым}, если его множество вершин содержится в
множестве вершин графа $G_{\tilde f_1}$, а внутренность каждого его
ребра содержится в $Z_\ell$; $Z_\ell$-допустимый однореберный
ориентированный граф $\gamma_0$ с параметризацией
$\gamma_0\:[0;1]\to\overline{Z_\ell}$ назовем {\it
$(Z_\ell,u)$-минимальным}, если
$$
u([\gamma_0])=\min\left\{u([\gamma])\ \left| \ \gamma\in
Z_\ell^\adm,\ \gamma|_{[0;1/2]}=\gamma_0|_{[0;1/2]},\
u([\gamma])>0\right.\right\},
$$
где через $Z_\ell^\adm$ обозначено множество $Z_\ell$-допустимых
однореберных графов. Из включений $u\in U_f^\infty \subset U_{\tilde
f_1}^\infty$ и определения подмножества $U_{\tilde f_1}^\infty\subset
H^1_f$ следует, что указанный минимум достигается для любого
$Z_\ell$-допустимого однореберного графа $\gamma_0$, причем ровно на
одном $Z_\ell$-допустимом однореберном графе с точностью до гомотопии
в классе $Z_\ell$-допустимых графов. Более того, имеется
$Z_\ell$-допустимый ориентированный граф
$G_\ell\subset\overline{Z_\ell}$, каждое ребро которого является
$(Z_\ell,u)$-минимальным и который содержит в качестве ровно одного
из своих ребер любой $(Z_\ell,u)$-минимальный однореберный граф с
точностью до гомотопии в классе $Z_\ell$-допустимых графов. Причем
граф $G_\ell$ с указанным свойством единствен с точностью до
гомотопии в классе $Z_\ell$-допустимых графов. Обозначим через
$G\subset M$ граф, полученный из графа $G_{\tilde f_1}$ заменой
подграфа $\d Z_\ell$ соответствующим графом $G_\ell$ для каждой
компоненты связности
$Z_\ell\subset M\setminus G_{\tilde f_1}$, такой что
$Z_\ell\subset\tilde f_1^{-1}([c^-;c^+])$. Нетрудно показывается, что
$G_f\in(\Diff^0(M,\N_f))(G)$.

В общем случае имеются такие последовательности граней
$\tau',\tau'',\dots,\tau^{(j-1)}$ грани $\tau=:\tau^{(j)}$ и
соответствующих отношений частичного порядка $\hat J=:J',J''$, $\dots$,
$J^{(j-1)}$, $J^{(j)}:=J$, что $\tau^{(i-1)}=\tau_{J^{(i-1)}}$
является гипергранью грани $\tau^{(i)}=\tau_{J^{(i)}}$, $2\le i\le
j$. Пусть $\tilde f_j:=f$ -- невозмущенная функция, и $\tilde f_i$ --
возмущенная функция, которой отвечает грань $\tau^{(i)}$ по
правилу~(\ref {eq:incid}), причем $\N_{\tilde f_i}=\N_f$, $1\le i\le
j-1$. Определим индуктивно по графу $G':=G_{\tilde f_1}=G_{\tilde
f^*}\subset M$ графы
$G^{(i+1)}:=G_{M,\N_{f,0,2},G^{(i)},u,J^{(i)},J^{(i+1)}}$ при
$i=1,2,\dots,j-1$. Положим $G_{M,\N_{f,0,2},G_{\tilde f_1},u,\hat
J,J}:=G^{(j)}$. Из доказанного выше следует (по индукции), что
$G_{\tilde f_{i+1}}\in(\Diff^0(M,\N_f))(G^{(i+1)})$ при
$i=1,2,\dots,j-1$. В частности, $G_f=G_{\tilde
f_j}\in(\Diff^0(M,\N_f))(G^{(j)})$. Лемма~\ref {lem:tilde} доказана.
\end{proof}

Из леммы~\ref {lem:tilde} следует ввиду~(\ref {eq:u}) и~(\ref
{eq:tilde}), что
$$
 G_f\in(\Diff^0(M,\N_f))(G_{M,\N_{f,0,2},G_{\tilde f^*},u_1,\hat J,J}),
 $$
 $$
 G_f\in(\Diff^0(M,\N_f))(G_{M,\N_{f,0,2},G_{\tilde f^*_1},u_2,J(\barc(\tilde f_1^*)),J}),
$$
 $$
 G_{fh_1}\in(\Diff^0(M,\N_{fh_1}))(G_{M,h_1^{-1}(\N_{f,0,2}),h_1^{-1}(G_{\tilde f^*_1}),h_1^*(u_2),J(\barc(\tilde f^*_1h_1)),J(\barc(fh_1))})
 $$
 $$
 =(\Diff^0(M,\N_f))(G_{M,\N_{f,0,2},G_{\tilde f^*},u_1,\hat J,J})
 =(\Diff^0(M,\N_f))(G_f),
 $$
т.е.\ следует требуемое включение
$G_{fh_1}\in(\Diff^0(M,\N_f))(G_f)$. Таким образом,
$h_1\in(\stab_{\D^0}f)(\Diff^0(M,\N_f))$, откуда
$u_1\in\widetilde\Gamma_f^*(u_2)$, а значит, отображение $\chi_{f,g}$
является вложением.

\section{Построение комплекса $\KK$ оснащенных функций Морса} \label {sec:KK}

В данном параграфе строится косой цилиндрически-полиэдральный
комплекс $\KK$ (см.\ определение~\ref {def:pol}(B)), удовлетворяющий
условиям теоремы~\ref {thm:KP4}.
Мы получим комплекс $\KK$ из стандартных косых цилиндрических ручек,
описанных в \S\ref{sec:KK0}, путем приклеивания друг к другу по
построенным там отображениям инцидентности. Конструкция, приведенная
в данном параграфе, является обобщением конструкции
из~\cite[\S5]{BaK}, где были построены (при дополнительном
ограничении $p=p^*$, $q=q^*$, $r=r^*$, т.е.\ когда все критические
точки фиксированы) более простые комплексы $\K$ и $K$ (комплексы
функций Морса), являющиеся строго полиэдральными комплексами (см.\
определение~\ref {def:pol},(C)).

Пусть $f_*\in F^1$ -- базисная функция Морса (см.\ определение~\ref
{def:F}(A)). В каждом классе изотопности $[f]_\isot$ рассмотрим
отмеченную функцию $f\in[f]_\isot$ этого класса с множествами
критических точек $\N_{f,\lam}=\N_{f_*,\lam}$, $\lam=0,1,2$ (см.\
\S\ref {sec:KK0}, шаг 2), тогда $H^0_f=H^0_{f_*}$, $H^1_f=H^1_{f_*}$.
Рассмотрим топологическое пространство
$$
 (F^1/\sim_\isot)^\diskr \ \times \P^{q-1}_{f_*}\times H^1_{f_*}
 \ \approx \
 (F^1/\sim_\isot)^\diskr \ \times \P^{q-1}\times\RR^{2q},
$$
где $(F^1/\sim_\isot)^\diskr:=F^1/\sim_\isot$ с дискретной
топологией, а евклидов многогранник $\P^{q-1}_{f_*}\approx\P^{q-1}$ и
векторное пространство $H^1_{f_*}\approx\RR^{2q}$ определены как в
\S\ref {sec:KK0}, шаги 1, 4 и (\ref {eq:H1f}). Рассмотрим в этом
топологическом пространстве подпространство
 \begin{equation} \label {eq:XX}
 \XX:=\bigcup\limits_{[f]_\isot\in F^1/\sim_\isot}\{[f]_\isot\}\times D_{[f]_\isot}
 \times U_{[f]_\isot}
 \end{equation}
с индуцированной топологией (см.\ \S\ref{sec:KK0}, шаг 2 и~(\ref
{eq:U})).

\begin{Def} [комплекс $\KK$ оснащенных функций Морса] \label {def:KK}
Рассмотрим факторпространство $\KK:=(\XX/\sim)/\sim_\glue$ с
фактортопологией, где отношения эквивалентности $\sim$, $\sim_\glue$
на пространствах $\XX$, $\YY:=\XX/\sim$ порождены следующими
отношениями соответственно:

(отношение $\sim$ на $\XX$; стандартная косая цилиндрическая ручка
$\DD_{[f]_\isot}^\st$) для каждого класса изотопности $[f]_\isot$
рассмотрим проекцию $D_{[f]_\isot} \times U_{[f]_\isot}
\to(D_{[f]_\isot} \times
U_{[f]_\isot})/\widetilde\Gamma_{[f]_\isot}=\DD^\st_{[f]_\isot}$ на
стандартную косую цилиндрическую ручку $\DD^\st_{[f]_\isot}$ индекса
$k=\dim D_{[f]_\isot}$ (см.\ \S\ref {sec:KK0}, шаг 9), рассмотрим
проекцию $\{[f]_\isot\}\times D_{[f]_\isot}\times
U_{[f]_\isot}\to\{[f]_\isot\}\times\DD^\st_{[f]_\isot}=:\uups_{[f]_\isot}$,
являющуюся прямым произведением отображения
$\{[f]_\isot\}\to\{[f]_\isot\}$ и этой проекции, и назовем точки
множества $\{[f]_\isot\}\times D_{[f]_\isot} \times U_{[f]_\isot}$
{\it $\sim$-эквивалентными}, если их образы в $\uups_{[f]_\isot}$ при
последней проекции совпадают;
тогда $\YY:=\XX/\sim=\bigcup_{[f]_\isot}\uups_{[f]_\isot}$;

(отношение $\sim_\glue$ на $\YY$; отображения инцидентности) для каждой пары примыкающих классов $[f]_\isot\prec[g]_\isot$
рассмотрим вложение, называемое {\it отображением
инцидентности}: \\
$\chi_{[f]_\isot,[g]_\isot}\:\d_{[g]_\isot}\DD^\st_{[f]_\isot}\hookrightarrow\DD^\st_{[g]_\isot}$
(см.\ \S\ref {sec:KK0}, шаги 3 и 10), где
$\d_{[g]_\isot}\DD^\st_{[f]_\isot}$ содержится в подошве
$\d\DD^\st_{[f]_\isot}$ стандартной косой цилиндрической ручки
$\DD^\st_{[f]_\isot}$ и является объединением ее попарно
непересекающихся косых граней
$\d_{\tau'}\DD^\st_{[f]_\isot}\subset\d\DD^\st_{[f]_\isot}$,
см.~(\ref {eq:delta:DD}), (\ref {eq:delta:g:DD}); рассмотрим
индуцированное вложение
$\d_{[g]_\isot}\uups_{[f]_\isot}:=\{[f]_\isot\}\times(\d_{[g]_\isot}\DD^\st_{[f]_\isot})\hookrightarrow\uups_{[g]_\isot}$
(которое тоже обозначим через $\chi_{[f]_\isot,[g]_\isot}$); назовем
любую точку множества $\d_{[g]_\isot}\uups_{[f]_\isot}\subset\YY$ и
ее образ в $\uups_{[g]_\isot}\subset\YY$ при данном вложении {\it
$\sim_\glue$-эквивалентными}.
\end{Def}

Пусть $p_Y\:\YY\to\KK$ -- каноническая проекция. Рассмотрим
подмножество $\o\YY:=\bigcup\limits_{ [f]_\isot}
\o\uups_{[f]_\isot}\subset\YY$, где
$\o\uups_{[f]_\isot}:=\{[f]_\isot\}\times\o\DD_{[f]_\isot}^\st$ (см.\
определение~\ref {def:thick:cylinder},(C)). Обозначим
$\DD_{[f]_\isot}:=p_Y(\uups_{[f]_\isot})$,
$\o\DD_{[f]_\isot}:=p_Y(\o\uups_{[f]_\isot})$,
$\d\DD_{[f]_\isot}:=p_Y(\d\uups_{[f]_\isot})$,
$\d_{\tau'}\uups_{[f]_\isot}:=\{[f]_\isot\}\times(\d_{\tau'}\DD^\st_{[f]_\isot})$,
и через $\o\d_{[g]_\isot}\uups_{[f]_\isot}$ обозначим объединение
открытых (косых) граней
$\{[f]_\isot\}\times\left(((\o\tau')\times\SS_{[f]_\isot})/\Gamma_{[f]_\isot}\right)$
стандартной (косой) цилиндрической ручки $\uups_{[f]_\isot}$, таких
что $\tau'\prec D_{[f]_\isot}$ и $\delta_{\tau'}[f]_\isot=[g]_\isot$.

\begin{Thm} \label{thm:KK}
Пространство $\KK=\YY/\sim_\glue$ обладает структурой косого
цилиндрически-поли\-эд\-раль\-ного комплекса ранга $q-1$ с косыми
цилиндрическими ручками
$\DD_{[f]_\isot}=p_Y(\uups_{[f]_\isot})\subset\KK$, $[f]_\isot\in
F^1/\sim_\isot$. При этом для любого класса изотопности $[f]_\isot$
отображение
$\varphi_{[f]_\isot}:=p_Y|_{\uups_{[f]_\isot}}\:\uups_{[f]_\isot}\to\KK$
является характеристическим
отображением ручки $\DD_{[f]_\isot}$ (откуда отображение
$p_Y|_{\o\YY}\:\o\YY\to\KK$ биективно), и выполнено
 \begin{equation} \label {eq:prec'}
 \d\DD_{[f]_\isot} \subset \bigcup\limits_{[g]_\isot\succ[f]_\isot}\DD_{[g]_\isot}, \quad
 \DD_{[f]_\isot} \cap \o\DD_{[g]_\isot} = \varphi_{[f]_\isot}(\o\d_{[g]_\isot}\uups_{[f]_\isot}),
 \end{equation}
 \begin{equation} \label {eq:prec}
 \varphi_{[g]_\isot}\circ\chi_{[f]_\isot,[g]_\isot}=\varphi_{[f]_\isot}|_{\d_{[g]_\isot}\uups_{[f]_\isot}}
 \qquad \mbox{для любых} \quad [f]_\isot\prec[g]_\isot.
 \end{equation}
Дискретная группа $\D/\D^0$ действует на $\KK$ автоморфизмами косого
цилинд\-ри\-чес\-ки-по\-ли\-э\-драль\-ного комплекса.
\end{Thm}

\begin{proof}{} {\it Шаг 1.} При любом $k\in\ZZ$ рассмотрим подмножества
$$
   \YY^{(k)}:=\bigcup\limits_{\dim D_{[f]_\isot}\le k}  \uups_{[f]_\isot}\ \subset \YY = \XX/\sim, \quad
 \o\YY^{(k)}:=\bigcup\limits_{\dim D_{[f]_\isot}\le k}\o\uups_{[f]_\isot}\ \subset \YY^{(k)}
$$
с индуцированной топологией, и множество
$\KK^{(k)}:=\YY^{(k)}/\sim_\glue$ с фактортопологией. Докажем лемму
(индукцией по $k$) для пространств $\YY^{(k)}$, $\KK^{(k)}$,
$\o\YY^{(k)}\subset\YY^{(k)}$ и проекции
$p_{Y,k}:=p_Y|_{\YY^{(k)}}\:\YY^{(k)}\to\KK^{(k)}$. При $k<0$
доказывать нечего, так как $\KK^{(k)}=\emptyset$.

Пусть $k\ge1$, и доказываемое утверждение верно для $\YY^{(k-1)}$,
$\KK^{(k-1)}$. Покажем, что для каждого $[f]_\isot$, такого что
$\ind\uups_{[f]_\isot}=k$, имеется (``приклеивающее'') отображение
$\varphi'_{[f]_\isot}\:\d\uups_{[f]_\isot} \to\KK^{(k-1)}$ подошвы
$\d\uups_{[f]_\isot}$ косой цилиндрической ручки $\uups_{[f]_\isot}$,
такое что
$\varphi'_{[f]_\isot}|_{\d_{[g]_\isot}\uups_{[f]_\isot}}=p_{Y,k-1}\circ\chi_{[f]_\isot,[g]_\isot}$
для любого $[g]_\isot\succ[f]_\isot$ (а потому
$\ind\uups_{[g]_\isot}<k$), см.~(\ref {eq:delta:g:DD}) и~(\ref
{eq:delta:DD}). Действительно, это отображение однозначно,
так как для любых $[g_1]_\isot\succ[g]_\isot\succ[f]_\isot$ и
$\tau''\prec D_{[f]_\isot}$, таких что
$\delta_{\tau''}[f]_\isot=[g_1]_\isot$, выполнено
$p_{Y,k-1}\circ\chi_{[f]_\isot,[g_1]_\isot}|_{\d_{\tau''}\uups_f}
=p_{Y,k-1}\circ\chi_{[g]_\isot,[g_1]_\isot}\circ\chi_{[f]_\isot,[g]_\isot}|_{\d_{\tau''}\uups_f}
=p_{Y,k-1}\circ\chi_{[f]_\isot,[g]_\isot}|_{\d_{\tau''}\uups_f}$
ввиду (\ref{eq:obstruction:no}) и (\ref {eq:prec}) для $\YY^{(k-1)}$,
$\KK^{(k-1)}$. Отображение $\varphi'_{[f]_\isot}$ непрерывно, так как
его область определения является конечным (ввиду (\ref {eq:incid}))
объединением замкнутых подмножеств $\d_{[g]_\isot}\uups_{[f]_\isot}$,
таких что $[g]_\isot\succ[f]_\isot$, и его ограничение на такое
подмножество есть композиция
$p_{Y,k-1}\circ\chi_{[f]_\isot,[g]_\isot}$ непрерывных
отображений. Отображение $\varphi'_{[f]_\isot}$ инъективно, так как
его ограничение
$\varphi'_{[f]_\isot}|_{\o\d_{[g]_\isot}\uups_{[f]_\isot}}
=p_{Y,k-1}\circ\chi_{[f]_\isot,[g]_\isot}|_{\o\d_{[g]_\isot}\uups_{[f]_\isot}}$
на объединение открытых граней, отвечающих классу
$[g]_\isot\succ[f]_\isot$, является композицией вложений
$\chi_{[f]_\isot,[g]_\isot}|_{\o\d_{[g]_\isot}\uups_{[f]_\isot}}$ и
$p_{Y,k-1}|_{\o\uups_{[g]_\isot}}$ (ввиду \S\ref {sec:KK0}, шаг 10, и
биективности $p_{Y,k-1}|_{\o\YY^{(k-1)}}$), а образы таких
ограничений содержатся в открытых ручках
$p_{Y,k-1}\circ\chi_{[f]_\isot,[g]_\isot}(\o\d_{[g]_\isot}\uups_{[f]_\isot})\subset
p_{Y,k-1}(\o\uups_{[g]_\isot})$, которые попарно не пересекаются (для
разных классов $[g]_\isot$) ввиду биективности
$p_{Y,k-1}|_{\o\YY^{(k-1)}}$.

Отображение топологических пространств назовем {\it хорошим}, если
оно переводит любое замкнутое подмножество в замкнутое подмножество.
Отображение $\varphi'_{[f]_\isot}$ является хорошим, так как его
ограничение на каждую косую грань
$\d_{\tau'}\uups_{[f]_\isot}\subset\d_{[g]_\isot}\uups_{[f]_\isot}$
(для $\tau'\prec D_{[f]_\isot}$, $\delta_{\tau'}[f]_\isot=[g]_\isot$)
есть композиция
$p_{Y,k-1}|_{\uups_{[g]_\isot}}\circ\chi_{[f]_\isot,\tau'}$
(автоматически хорошего, см.\ определение~\ref {def:thick:cylinder})
мономорфизма $\chi_{[f]_\isot,\tau'}$ косых цилиндрических ручек
(см.\ \S\ref {sec:KK0}, шаг 10) и (автоматически хорошего)
характеристического отображения $p_{Y,k-1}|_{\uups_{[g]_\isot}}$
косой ручки $p_{Y,k-1}(\uups_{[g]_\isot})\subset\KK^{(k-1)}$ (по
предположению индукции), а каждая косая грань замкнута и их конечное
число. Отсюда следует, что $\KK^{(k)}$ гомеоморфно пространству (с
фактортопологией), полученному из $\KK^{(k-1)}$ приклеиванием косых
цилиндрических ручек $\uups_{[f]_\isot}$ индекса $k$ при помощи
инъективных непрерывных хороших отображений
$\varphi'_{[f]_\isot}\:\d\uups_{[f]_\isot}\to\KK^{(k-1)}$ их подошв,
а потому $\KK^{(k)}$ удовлетворяет условию~(w) из определения~\ref
{def:pol}(B). Поэтому $\KK^{(k)}$ обладает свойствами~(\ref
{eq:prec'}) и (\ref {eq:prec}), отображение
$p_{Y,k}|_{\uups_{[f]_\isot}}$ является инъективным, непрерывным и
хорошим, а потому задает структуру косой цилиндрической ручки на
$p_{Y,k}(\uups_{[f]_\isot})$ и является характеристическим
отображением этой ручки. Пространство $\KK^{(k)}$ с полученным
разбиением на косые цилиндрические ручки $p_{Y,k}(\uups_{[f]_\isot})$
удовлетворяет условию~(c) из определения~\ref {def:pol}(B), так как
ограничение характеристического отображения
$p_{Y,k}|_{\uups_{[f]_\isot}}$ любой ручки $\uups_{[f]_\isot}$
индекса $k$ на каждую свою косую грань $\d_{\tau'}\uups_{[f]_\isot}$
есть композиция
$p_{Y,k}|_{\d_{\tau'}\uups_{[f]_\isot}}=i_{k-1}\circ\varphi'_{[f]_\isot}|_{\d_{\tau'}\uups_{[f]_\isot}}
=i_{k-1}\circ p_{Y,k-1}\circ\chi_{[f]_\isot,\tau'}
=p_{Y,k}|_{\uups_{[g]_\isot}}\circ\chi_{[f]_\isot,\tau'}$
мономорфизма $\chi_{[f]_\isot,\tau'}$ стандартных косых
цилиндрических ручек и характеристического отображения
$p_{Y,k}|_{\uups_{[g]_\isot}}$ косой цилиндрической ручки
$p_{Y,k}(\uups_{[g]_\isot})$, где
$i_{k-1}\:\KK^{(k-1)}\hookrightarrow\KK^{(k)}$ -- отображение
включения. Значит, $\KK^{(k)}$ является косым
цилиндрически-полиэдральным комплексом ранга $k$, что завершает
доказательство индукционного перехода.

{\it Шаг 2.} Определим (естественное) правое действие
$\rho_Y\:\D/\D^0\times\YY\to\YY$ дискретной группы $\D/\D^0$ на
пространстве $\YY$ формулой
$$
 (h\D^0,[f]_\isot,\widetilde\Gamma_{[f]_\isot}(\barc,u))\mapsto
 ([fh]_\isot,\widetilde\Gamma_{[fh]_\isot}(h_{1;[f]_\isot,[fh]_\isot}^{*0}(\barc),h_{1;[f]_\isot,[fh]_\isot}^*(u)))
$$
при любых $h\in\D$, $([f]_\isot,\barc,u)\in\XX$ (см.\ \S\ref
{sec:KK0}, конец шага 9), где $h_{1;[f]_\isot,[fh]_\isot}\in
h\D^0\subset\D$ -- диффеоморфизм, переводящий линии уровня
отмеченной функции класса изотопности $[fh]_\isot$ в линии уровня
отмеченной функции класса изотопности $[f]_\isot$ с сохранением
направления роста. Это действие определено корректно, так как в
силу~\cite[лемма~1]{KP2} для любых $h_1,h_2\in\D$ выполнено
$$
h_{1;[f]_\isot,[fh_1h_2]_\isot}^{-1}h_{1;[f]_\isot,[fh_1]_\isot}h_{1;[fh_1]_\isot,[fh_1h_2]_\isot}\in(\stab_{\D^0}{g})(\Diff^0(M,\N_{g}))
$$
и действие групп $\stab_{g}\D^0$ и $\Diff^0(M,\N_{g})$ на косой
цилиндрической ручке $\DD^\st_{[g]_\isot}$ тривиально (см.\ (\ref
{eq:DD}), (\ref {eq:Gamma:f})), где $g$ -- отмеченная функция класса
изотопности $[fh_1h_2]_\isot$. Для любых $[f]_\isot\prec[g]_\isot$ (а
потому $[fh]_\isot\prec[gh]_\isot$ для любого $h\in\D$) определим
отображение инцидентности
$$
\chi_{[f],[g]}\:\bigcup\limits_{h\in\D}(\d_{[gh]_\isot}\uups_{[fh]_\isot})\to\bigcup\limits_{h\in\D}\uups_{[gh]_\isot}
$$
правилом
$\chi_{[f],[g]}|_{\d_{[gh]_\isot}\uups_{[fh]_\isot}}:=\chi_{[fh]_\isot,[gh]_\isot}\:\d_{[gh]_\isot}\uups_{[fh]_\isot}\to\uups_{[gh]_\isot}$
для любого $h\in\D$. Отображение $\chi_{[f],[g]}$ определено
корректно, так как для любых $[f]_\isot$ и $h_1\in\stab_{\D^0}f$
подмножества $\d_{[g]_\isot}\uups_{[f]_\isot}$,
$\d_{[gh_1]_\isot}\uups_{[f]_\isot}\subset\d\uups_{[f]_\isot}$ либо
совпадают (откуда $[g]_\isot=[gh_1]_\isot$), либо не пересекаются
ввиду леммы~\ref {lem:P} (см.\ \S\ref {sec:KK0}, шаг 10). Отображение
$\chi_{[f],[g]}$ является $(\D/\D^0)$-экви\-ва\-ри\-ант\-ным (т.е.\
$\chi_{[f],[g]}\circ\rho_Y(h,\cdot)=\rho_Y(h,\chi_{[f],[g]}(\cdot))$
для любого $h\in\D$), так как ввиду (\ref {eq:corr}) выполнено \\
$h^{-1}h_{[f]_\isot,\tau'}^{-1}hh_{[fh]_\isot,h^*_{1;[f]_\isot,[fh]_\isot}(\tau')}\in(\stab_{\D^0}{g_1})(\Diff^0(M,\N_{g_1}))$
и действия групп $\stab_{\D^0}{g_1}$ и $\Diff^0(M,\N_{g_1})$ на
стандартной косой цилиндрической ручке $\DD^\st_{[g_1]_\isot}$
тривиальны, где $\tau'\prec D_{[f]_\isot}$,
$\delta_{\tau'}[f]_\isot=[g]_\isot$, $g_1$ -- отмеченная функция
класса изотопности $[gh]_\isot$, $h\in\D$. Поэтому правое действие
$\rho_Y$ индуцирует правое действие $\rho^*\:\D/\D^0\times\KK\to\KK$
автоморфизмами косого цилиндрически-полиэдрального комплекса, такое
что $\rho^*(h\D^0,p_Y(y)):=p_Y\circ\rho_Y(h\D^0,y)$ для любых
$h\in\D$, $y\in\YY$.

В частности, действие $\rho^*(h\D^0,\cdot)$ любого элемента
$h\D^0\in\D/\D^0$ на комплексе $\KK$ индуцирует изоморфизм
$\DD_{[f]_\isot}\to\DD_{[fh]_\isot}$ косых цилиндрических ручек
$\DD_{[f]_\isot},\DD_{[fh]_\isot}\subset\KK$ для любой функции $f\in
F$. Поэтому эти ручки изоморфны одной и той же стандартной ручке
$(D_{f_0}\times\SS_{f_0})/\Gamma_{f_0}$, где $f_0$ -- отмеченная
функция класса эквивалентности $[f]=[fh]$.
Теорема \ref {thm:KK} полностью доказана.
\end{proof}


Аналогично определению \ref {def:KK} определим гладкое многообразие
$\MM$. А именно, рассмотрим в евклидовом пространстве
$H^0_f=\RR^{\N_{f,1}}\cong\RR^q$ открытый куб
$(-1;1)^{\N_{f,1}}\cong(-1;1)^q$, и для любой грани
$\tau=\tau_J\subset\P_f^{q-1}$ рассмотрим ее внутренность $\o\tau$ и
обозначим через $(\o\tau)^*$ множество таких 0-коцепей
$\barc'\in(-1;1)^{\N_{f,1}}\subset H^0_f\cong\RR^q$, что
$J(\barc')=J$ (см.\ \S\ref {sec:KK0}, шаг 1), а через $\tau^*$
обозначим множество таких 0-коцепей
$\barc'\in(-1;1)^{\N_{f,1}}\cong(-1;1)^q$, что $J(\barc')\preceq J$.
Назовем $(\o\tau)^*$ {\em стратом}, а $\tau^*$ -- {\em звездой} этого
страта (отвечающими грани $\tau$). Тогда звезда $\tau^*$ открыта в
$H^0_f$,
 \begin{equation} \label {eq:strat}
 (\o\tau)^*\subseteq\tau^*, \qquad  \tau\subset\tau^*,
 \end{equation}
страт $(\o\tau)^*$ есть выпуклое открытое подмножество
некоторого $|f(\N_{f,1})|$-мерного линейного подпространства в
$H^0_f$, грань $\tau$ есть выпуклый $(q-|f(\N_{f,1})|)$-мерный
многогранник, причем грань $\tau$ и
страт $(\o\tau)^*$ пересекаются трансверсально в барицентре грани
$\tau$, являющемся внутренней точкой каждого из них. Положим
$D^*_{[f]_\isot}:=\tau_{J(\barc(f))}^*$, $(\o
D_{[f]_\isot})^*:=(\o\tau_{J(\barc(f))})^*$. Рассмотрим в
пространстве $(F^1/\sim_\isot)^\diskr\times (-1;1)^{\N_{f,1}} \times
H^1_{f_*}$ подпространства
$$
 \XX^{\infty,*}:=\bigcup\limits_{[f]_\isot\in F^1/\sim_\isot} \{[f]_\isot\}\times D^*_{[f]_\isot}
 \times U_{[f]_\isot}^{\infty},
$$
$$
 \XX^{\infty,\circ*}:=\bigcup\limits_{[f]_\isot\in F^1/\sim_\isot} \{[f]_\isot\}\times (\o D_{[f]_\isot})^*
 \times U_{[f]_\isot}^{\infty},
$$
где $U_{[f]_\isot}^{\infty}$ определено как в~(\ref {eq:U0}). Имеем
включения $\XX\subset \XX^{\infty,*}$ и
$\XX^{\infty,\circ*}\subset\XX^{\infty,*}$.

Аналогично определению~\ref {def:KK} определим отношение
эквивалентности $\sim$ на $\XX^{\infty,*}$ с помощью покомпонентного
действия на $D^*_{[f]_\isot}\times U_{[f]_\isot}^{\infty}$ дискретной группы $\widetilde\Gamma_{[f]_\isot}$ (см.\ \S\ref
{sec:KK0}, шаг 9).
Определим ``окрестность'' стандартной ручки $\DD_{[f]_\isot}^\st$:
$$
 (\DD_{[f]_\isot}^{\infty,*})^\st:=(D^*_{[f]_\isot}
\times
U_{[f]_\isot}^{\infty})/\widetilde\Gamma_{[f]_\isot}\approx(D^*_{[f]_\isot}\times\SS_{[f]_\isot}^{\infty})/\Gamma_{[f]_\isot},
$$
$$
\uups_{[f]_\isot}^{\infty,*}:=\{[f]_\isot\}\times(\DD_{[f]_\isot}^{\infty,*})^\st\subset\YY^{\infty,*}:=\XX^{\infty,*}/\sim.
$$
Определим отношение эквивалентности $\sim_\glue$ на
$\YY^{\infty,*}:=\XX^{\infty,*}/\sim$ с помощью вложений \\
$\chi_{[f]_\isot,[g]_\isot}^{\infty,*}\:\d_{[g]_\isot}(\DD_{[f]_\isot}^{\infty,*})^\st\hookrightarrow(\DD_{[g]_\isot}^{\infty,*})^\st$,
определяемых теми же формулами, что и вложения \\
$\chi_{[f]_\isot,[g]_\isot}\:\d_{[g]_\isot}\DD^\st_{[f]_\isot}\hookrightarrow\DD^\st_{[g]_\isot}$
(см.~\S\ref {sec:KK0}, шаг 10). Рассмотрим пространства
$$
 \YY^{\infty,\circ*}:=\XX^{\infty,\circ*}/\sim, \qquad
 \YY^{\infty,*}:=\XX^{\infty,*}/\sim, \qquad
 \MM:=\YY^{\infty,*}/\sim_\glue
 $$
с фактортопологией, тогда $\KK\subset \MM$.
Пусть
 $$
p_X\:\XX^{\infty,*}\to\MM, \qquad p_Y\:\YY^{\infty,*}\to\MM
 $$
-- канонические проекции. Рассмотрим ``окрестность''
$\DD_{[f]_\isot}^{\infty,*}:=p_Y(\uups_{[f]_\isot}^{\infty,*})$ косой
цилиндрической ручки $\DD_{[f]_\isot}\subset\KK$ в $\MM$.

Так как $\XX^{\infty,*}$ является гладким открытым $3q$-мерным
многообразием с естественной плоской аффинной связностью, и каждая
группа $\widetilde\Gamma_{[f]_\isot}$ действует на нем
диффеоморфизмами, сохраняющими связность, то $\YY^{\infty,*}$ тоже
является гладким открытым $3q$-мерным многообразием с плоской
аффинной связностью.

\begin{Thm} \label {thm:MM1}
Пространство $\MM:=\YY^{\infty,*}/\sim_\glue$ обладает структурой
гладкого $3q$-мерного многообразия и естественной плоской аффинной
связностью, гладкой относительно этой структуры. При этом для каждого
класса изотопности $[f]_\isot$ выполнены следующие условия:
\\
(i) подмножество
$\DD_{[f]_\isot}^{\infty,*}=p_Y(\uups_{[f]_\isot}^{\infty,*})\subset\MM$
открыто, где $p_Y\:\YY^{\infty,*}\to\MM$ -- проекция,
\\
(ii) отображение
$p_Y|_{\uups_{[f]_\isot}^{\infty,*}}\:\uups_{[f]_\isot}^{\infty,*}\to
\MM$ является гладким регулярным вложением гладких $3q$-мерных
многообразий, сохраняющим аффинную связность.
\\
Отображение $p_Y|_{\YY^{\infty,\circ*}}\:\YY^{\infty,\circ*}\to\MM$
биективно. Дискретная группа $\D/\D^0$ и группа $\Diff^+[-1;1]$
действуют на $\MM$ справа и слева соответственно диффеоморфизмами,
сохраняющими аффинную связность и области
$\DD_{[f]_\isot}^{\infty,*}$.
\end{Thm}

\begin{proof}{} Аналогично доказательству теоремы~\ref {thm:KK}
доказательство проводится индукцией, а именно, доказывается
однозначность, непрерывность и инъективность соответствующих
приклеивающих отображений, а также биективность отображения
$p_Y|_{\YY^{\infty,\circ*}}\:\YY^{\infty,\circ*}\to\MM$.
Так как каждое приклеивающее отображение
$(\varphi_{[f]_\isot}^{\infty,*})'$ определено на открытом
подмножестве открытого подмножества
$\uups_{[f]_\isot}^{\infty,*}\subset\YY^{\infty,*}$ и является
непрерывным и инъективным отображением $3q$-мерных многообразий, то
$\DD_{[f]_\isot}^{\infty,*}$ открыто в $\MM$ и
$p_Y|_{\uups_{[f]_\isot}^{\infty,*}}\:\uups_{[f]_\isot}^{\infty,*}\to\DD_{[f]_\isot}^{\infty,*}\subset\MM$
является гомеоморфизмом (в индуцированной топологии). В частности,
каждая точка пространства $\MM$ обладает окрестностью, гомеоморфной
$\RR^{3q}$. Так как пространство $\MM$ хаусдорфово (см.\ ниже), а все
отображения $\chi_{[f]_\isot,[g]_\isot}^{\infty,*}$ (а потому и все
приклеивающие отображения) являются гладкими и сохраняют плоскую
аффинную связность, то $\MM$ является гладким $3q$-мерным
многообразием с естественной плоской аффинной связностью.

Осталось доказать, что пространство $\MM$ хаусдорфово. Рассмотрим
естественное непрерывное ``вычисляющее'' отображение
 $$
\eval^*\:\XX^{\infty,*}\to\RR^q/\Sigma_q, \qquad
([f]_\isot,\barc,u)\mapsto\Sigma_q\barc,
 $$
где группа перестановок $\Sigma_q$ действует на любом наборе
$\barc\in\RR^{\N_{f,1}}\cong\RR^q$ перестановками компонент. Легко
проверяется, что оно индуцирует (автоматически непрерывное)
отображение $ \MM\to\RR^q/\Sigma_q$, которое тоже будем обозначать
через $\eval^*$. Пусть $m_1,m_2\in\MM$. Если
$\eval^*(m_1)\ne\eval^*(m_2)$, то точки $m_1,m_2$ обладают
непересекающимися окрестностями ввиду хаусдорфовости $\RR^q/\Sigma_q$
и непрерывности $\eval^*$. Пусть $\eval^*(m_1)=\eval^*(m_2)$ и
$m_i\in\DD_{[f_i]_\isot}^{\infty,*} $, $i=1,2$. Если
$[f_1]_\isot\ne[f_2]_\isot$, то в силу леммы~\ref {lem:P} не
существует класса изотопности $[g]_\isot$, такого что
$[g]_\isot\succ[f_i]_\isot$, $i=1,2$, а потому образы приклеивающих
отображений $(\varphi_{[f_i]_\isot}^{\infty,*})'$, $i=1,2$, имеют
пустое пересечение, откуда окрестности $\DD_{[f_i]_\isot}^{\infty,*}$
точек $m_i$, $i=1,2$, в $\MM$ имеют пустое пересечение. Если
$[f_1]_\isot=[f_2]_\isot$ и $m_1\ne m_2$, то точки
$m_1,m_2\in\DD_{[f_1]_\isot}^{\infty,*}$ обладают непересекающимися
окрестностями ввиду открытости $\DD_{[f_1]_\isot}^{\infty,*}$,
хаусдорфовости стандартного пространства 
$(\DD_{[f_1]_\isot}^{\infty,*})^\st$ и гомеоморфности
$\DD_{[f_1]_\isot}^{\infty,*}\approx\uups_{[f_1]_\isot}^{\infty,*}\approx(\DD_{[f_1]_\isot}^{\infty,*})^\st$ (см.\ выше). Таким образом, пространство $\MM$ является хаусдорфовым.
Теорема \ref {thm:MM1} доказана.
\end{proof}

Теоремы~\ref {thm:KK} и \ref {thm:MM1} доказывают основную
теорему~\ref {thm:KP4}(A,B). Утверждение п.(C) теоремы~\ref {thm:KP4}
нетрудно выводится из (\ref {eq:3.23.5}) и явного описания (\ref
{eq:generH}) набора
образующих свободной абелевой группы $\D^0\cap\tildeGamma_f$ (см.\
\S\ref {sec:KK0}, шаг 7).




\begin{thebibliography}{99}

\bibitem{kp1} Е.А.\ Кудрявцева, Д.А.\ Пермяков. Оснащенные функции
Морса на поверхностях. Матем.\ Сб. {201}, No.\,4 (2010), 501-567.


\bibitem{Gromov86} М.\ Громов. {Дифференциальные Соотношения с Частными
Производными}. М.: Мир, 1990.

\bibitem{Igusa84} K.\ Igusa. {Higher singularities of smooth functions are
unnecessary}. Ann.\ Math. {119} (1984), 1-58.

\bibitem{Vassil89} В.А.\ Васильев. {Топология пространств функций без сложных
особенностей}. Функ.\ Ан.\ Прил. {23}(4) (1989), 24-36.

\bibitem{Franks} J.\ Franks. {Homology and Dynamical Systems}. CBMS
Regional Conf.\ {49} (1982), Amer.\ Math.\ Soc.

\bibitem{Hatcher75} A.\ Hatcher. {Higher simple homotopy theory}.
Annals of Math. (2) {102} (1975), 101-137.

\bibitem{ChL09} A.\ Chenciner, F.\ Laudenbach. {Morse 2-jet space and
$h$-principle}. arXiv:0902.3692v1 [math.GT] 23 Feb 2009

\bibitem{BaK} Е.А.\ Кудрявцева. {Связные компоненты пространств функций Морса с
фиксированными критическими точками}. Вестн.\ Моск.\ Ун-та. Сер.~1,
Математика. Механика, в печати (2011). arXiv:1007.4398.

\bibitem{FI}
B.\ Farb, N.V.\ Ivanov. {The Torelli geometry and its applications.
Research announcement}. arXiv:math.GT/0311123 v1.

\bibitem{H} W.J.\ Harvey. {Geometric structure of surface mapping class groups}.
In Homological group theory (Proc. Sympos., Durham, 1977), 255-269.
Cambridge: Cambridge Univ. Press, 1979.

\bibitem{Ha} J.\ Harer. {The virtual cohomological dimension of the mapping class
group of an orientable surface}. Invent. Math. 84(1) (1986), 157-176.

\bibitem{I1} N.\ Ivanov. {Automorphisms of complexes of curves and of
Teichm\"uller spaces}. In: Progress in knot theory and related
topics, 113-120. Paris: Hermann, 1997.

\bibitem{HT} A.\
Hatcher, W.\ Thurston. {A presentation for the mapping class group of
a closed orientable surface}. Topology 19(3) (1980), 221-237.

\bibitem{Mar}
D.\ Margalit. {The automorphism group of the complex of pants
decompositions}. arXiv:math/0201319v1 [math.GT].

\bibitem{MP} S.V.\ Matveev, M.\ Polyak.
{Cubic complexes and finite types invariants}. Geom.\ Topol.\ Monogr.
{4} (2002), 215-233. arXiv: math.GT/0204085, 2002.

\bibitem{Kmsb} Е.А.\ Кудрявцева. {Реализация гладких функций на поверхностях в виде функций высоты}. Матем.\ Сборник
{190} (1999), No.\,3, 29-88.

\bibitem{SH} В.В.\ Шарко. {Функции на поверхностях}, I. В книге: Труды
Матем.\ Инст.\ Укр.\ НАН ``Некоторые проблемы современной
математики'', ред.\ В.В.Шарко, {25}, Киев, Наукова Думка, 1998.
С.~408-434.

\bibitem{Max2005} S.I.\ Maksymenko. {Path-components of Morse mappings spaces of
surfaces}. Comment.\ Math.\ Helv. {80} (2005), 655-690.

\bibitem{Bu} Ю.М.\ Бурман. {Теория Морса для функций двух переменных
без критических точек}. Функц.\ Дифф.\ Ур. {3}(1-2) (1995), 31-31.

\bibitem{Bu2} Yu.M.\ Burman. Triangulations of surfaces with boundary and the homotopy principle
for functions without critical points. Annals of Global Analysis and
Geometry {17}(3) (1999), 221-238.

\bibitem{Kul} E.V.\ Kulinich. {On topologically equivalent Morse functions on
surfaces}. Methods of Funct.\ Anal.\ Topology {4} (1) (1998), 59-64.

\bibitem{Max} S.I.\ Maksymenko. {Stabilizers and orbits of Morse
functions}. arXiv:math.GT/0310067 v5 14 Aug 2006.

\bibitem{F86} А.Т.\ Фоменко. {Теория Морса интегрируемых гамильтоновых систем}.
ДАН СССР 287, No.\,5 (1986), 1071-1075.

\bibitem{FZ0} А.Т.\ Фоменко, Х.~Цишанг. {Топологический инвариант и
критерий эквивалентности интегрируемых гамильтоновых систем с двумя
степенями свободы}. Изв.\ АН СССР {54}, No.\,3 (1990), 546-575.

\bibitem{BF0} А.В.\ Болсинов, А.Т.\ Фоменко. {Траекторная
эквивалентность интегрируемых гамильтоновых систем с двумя степенями
свободы. Теорема классификации}. I: Матем.\ Сб. {185}, No.\,4 (1994),
27-89; II: Матем.\ Сб. {185}, No.\,5 (1994), 27-28.

\bibitem{BFbook} А.В.\ Болсинов, А.Т.\ Фоменко. {Введение в топологию
интегрируемых гамильтонорвых систем}. М.: Наука, 1997.

\bibitem{K} Е.А.\ Кудрявцева. {Устойчивые топологические и
гладкие инварианты сопряженности гамильтоновых систем на
поверхностях}. В книге: Топологические методы в теории гамильтоновых
систем. Ред.\ А.Т.\ Фоменко и А.В.\ Болсинов. М.: Факториал, 1998.
C.~147-202.

\bibitem{K:stab} Кудрявцева Е.А., Устойчивые инварианты сопряженности
гамильтоновых систем на двумерных поверхностях. Докл. Акад. Наук 361,
N.3 (1998), 314-317.

\bibitem{BrK:stab} {Brailov, Yu.\,A. and Kudryavtseva, E.\,A.}, Stable
topological nonconjugacy of Hamiltonian systems on two-dimensional
surfaces. Vestnik Moskov. Univ. Ser. I Mat.\ Mekh. {\bf No.\,2}
(1999), 20-27, 72 (in Russian).

\bibitem{A89} В.И.\ Арнольд. {Пространства функций с умеренными
особенностями}. Функц.\ Анал.\ Прил. {23}(3) (1989), 1-10.


\bibitem {KP2} Е.А.\ Кудрявцева.
{Равномерная лемма Морса и критерий изотопности функций Морса на
поверхностях}. Вестн.\ Моск.\ Ун-та. Сер.~1, Математика. Механика,
No.\,4 (2009), 13-22.

\bibitem{EE} C.J.\ Earle, J.\ Eells, Jr. {The diffeomorphism group of
a compact Riemann surface}. Bull.\ Amer.\ Math.\ Soc. {73}, no.~4
(1967), 557-559.

\bibitem{EE0} C.J.\ Earle, J.\ Eells, Jr. {A fibre bundle description
of Teichm\"uller theory}. J.\ Diff.\ Geometry {3} (1969), 19-43.

\bibitem{S} S.\ Smale. {Diffeomorphisms of the 2-sphere}. Proc.\ Amer.\
Math.\ Soc. {10} (1959), 621-626.

\bibitem{BrHae} {\it Bridson M.R., Haefliger A.}, Metric spaces of
non-positive curvature // Berlin, Heidelberg, N.Y., Barcelona, Hong
Kong, London, Milan, Paris, Singapore, Tokyo: Springer, 1999.

\bibitem{FF} А.Т.\ Фоменко, Д.Б.~Фукс. {Курс гомотопической топологии}.
М.: Наука, 1989.

\bibitem{P} {A.\ Postnikov.} Permutohedra, associahedra, and beyond.
arXiv:math/0507163v1 [math.CO] 7 Jul 2005.

\bibitem{Dehn} M.\ Dehn. {Die Gruppe der Abbildungsklassen (Das arithmetische Feld auf
Fl\"achen)}. Acta mathematica {69} (1938), 135-206.

\bibitem{BLMC} J.S. Birman, A. Lubotzky, J. McCarthy. Abelian and solvable
subgroups of the mapping class group. Duke Math.\ J. 50, No.4 (1983),
1107-1120.

\bibitem{Pdiplom} {Д.А.\ Пермяков.} Линейная независимость скручиваний Дэна. Дипломная работа.
http://dfgm.math.msu.su/files/0students/2009-dip-permyakov.pdf

\bibitem{Kronrod} А.\ Кронрод. {О функциях двух переменных}. Успехи Матем.\ Наук {5} (1950), No.\,1, 24-134.

\end{thebibliography}
\end{document}